\documentclass[preprint,12pt]{elsarticle}

\newtheorem{theorem}{Theorem}[section]
\newtheorem{lemma}{Lemma}[section]
\newtheorem{corollary}{Corollary}[section]
\newtheorem{proposition}{Proposition}[section]

\newtheorem{remark}{Remark}[section]
\newcommand{\ignore}[1]{}{}

\usepackage{amssymb}
\usepackage{amsmath}

\usepackage{color}
\usepackage{soul}
\usepackage[dvipsnames]{xcolor}

\usepackage[colorlinks=true, linkcolor=blue, citecolor=blue]{hyperref}

\def\1{{{\mbox{${\rm{1\negthinspace\negthinspace I}}$}}}}

\newcommand\beq{\begin{equation}}
\newcommand\eeq{\end{equation}}

\usepackage{ifthen}
\usepackage{xkeyval}
\usepackage{todonotes}
\begin{document}

\begin{frontmatter}

\title{Self-normalized Cram\'{e}r type moderate deviations for stationary sequences and applications}
\author{Xiequan Fan, \ \   Ion Grama,\  \  Quansheng Liu \ \   and\ \    Qi-Man Shao }
\address{Center for Applied Mathematics,
Tianjin University, Tianjin, China. }
\address{Universit\'{e} de Bretagne-Sud, LMBA, UMR CNRS 6205,
 Campus de Tohannic,\\ 56017 Vannes, France.}
\address{Department of Statistics and Data Science, Southern University of Science and Technology,\\
 Shenzhen  518000, China.}

\begin{abstract}
Let $(X _i)_{i\geq1}$ be a stationary sequence.
Denote $m=\lfloor n^\alpha \rfloor, 0< \alpha < 1,$ and $ k=\lfloor n/m \rfloor,$
 where  $\lfloor a \rfloor$ stands for the integer part of $a.$
Set $S_{j}^\circ     = \sum_{i=1}^m X_{m(j-1)+i}, 1\leq j \leq k,$
and $  (V_k^\circ)^2 =  \sum_{j=1}^k (S_{j}^\circ)^2 .$
We prove a Cram\'er type moderate deviation expansion  for $\mathbb{P}(  \sum_{j=1}^k S_{j}^\circ /V_k^\circ   \geq x)$ as $n\to \infty.$
Applications to mixing type sequences,   contracting Markov chains, expanding maps  and confidence intervals  are discussed.
\end{abstract}

\begin{keyword} moderate deviations;  stationary processes; Cram\'{e}r  moderate deviations
\vspace{0.3cm}
\MSC primary 60F10; 60G10; secondary  60E15
\end{keyword}

\end{frontmatter}




\section{Introduction}
 Let  $(X_i)_{i\geq 1}$ be a sequence of independent and identically distributed (i.i.d.) centered real random variables,
 that is $ \mathbb{E} X_{1}=0.$
Denote $S_n=\sum_{i=1}^{n}X_{i}$ the partial sums of $(X_i)_{i\geq 1}$ and $\sigma^2=\mathbb{E}X_{1}^2$ the variance of $X_1$.  Cram\'{e}r \cite{Cramer38} has established the following  asymptotic  moderate deviation expansion for the standardized sums:  if $\mathbb{E}\exp\{ c_{0}|X_{1}|\}<\infty$ for some constant $c_{0}>0$, termed Cram\'er's condition,  then for all $0\leq x =  o(n^{1/2}), $
\begin{equation}
\Bigg| \ln \frac {\mathbb{P}(S_n\geq x\sigma\sqrt{n})} {1-\Phi(x)} \Bigg|  = O(1)  \frac{(1+x)^3}{\sqrt{n}}  \ \ \mbox{as} \ \ n \rightarrow \infty,
\label{cramer1}
\end{equation}
where $\Phi(x)=\frac{1}{\sqrt{2\pi}}\int_{-\infty}^{x}\exp\{-t^2/2\}dt$ is the standard normal distribution.
Inequality (\ref{cramer1}) implies that
\begin{equation}\label{Cramer00}
\frac{\mathbb{P}(S_n\geq x\sigma\sqrt{n})}{1-\Phi \left( x\right)} = 1+ o(1)
\end{equation}
uniformly in the normal range $0\leq x =o( n^{ 1/6 } ) .$  Notice that Cram\'er's condition is  sufficient but not necessary for
(\ref{Cramer00}) to hold.
Indeed,   Linnik \cite{L61}   proved that for $\alpha \in (0, \frac16],$ formula (\ref{Cramer00})
holds uniformly for $0\leq x =o(n^{\alpha})$  as $n\rightarrow \infty$ if and only if $\mathbb{E}e^{ |X_1|^{4\alpha/(2\alpha+1)} } < \infty.$
Following the seminal work of Cram\'{e}r, various moderate deviation expansions for standardized  sums  have been obtained by many authors, see, for instance, Petrov \cite{Pe54}, Saulis and Statulevi\v{c}ius \cite{SS78} and \cite{F17}.
See also  Ra\v{c}kauskas \cite{Rackauskas95,Rackauskas97}, Grama \cite{G97},  Grama and Haeusler \cite{GH00} and \cite{FGL13} for martingales,
and  Wu and Zhao \cite{WZ08}  and Cuny and Merlev\`{e}de \cite{CM14} for stationary processes.

For establishing  moderate deviation expansions of type (\ref{Cramer00}) with a range $0\leq x =o( n^{ \alpha} ),$ $\alpha>0$, Linnik's condition   is necessary. However, Linnik's condition  becomes too restrictive if we only have finite moments of order $2+\rho, \rho \in (0, 1]$.
Although we still can
 establish (\ref{Cramer00}) via (non-uniform) Berry-Esseen estimations (see Bikelis \cite{B66}), the range cannot be wider than $0\leq x =O( \sqrt{\ln n}  ) ,$
 which is much more narrow than $0\leq x =o( n^{ \alpha} )$.
 To overcome this limitation, 
 instead of considering the   standardized sums,  one may consider the self-normalized sums, defined as follows:
  $$W_n= S_n/V_n,\ \ \ \ \ \ \ \textrm{where}\   \
 V_n^2= \sum_{i=1}^n X_i^2.$$
One of the motivations to consider self-normalized  sums is
 due to Student's $t$-statistic: 
\[
T_n = \sqrt{n} \, \overline{X}_n / \widehat{\sigma},
\]
where $$\overline{X}_n = \frac{S_n}{n}  \ \ \ \ \ \ \ \ \textrm{and}\ \  \ \ \ \ \ \ \widehat{\sigma}^2 = \sum_{i=1}^n  \frac{(X_i - \overline{X}_n )^2}{ n-1}  .$$
The Student's $t$-statistic $T_n$ is linked to the the self-normalized sum $W_n$ by the following formula: for all $x\geq0,$
\[
\mathbb{P}\Big( T_n  \geq x \Big) = \mathbb{P}\bigg(  W_n  \geq x \Big(\frac{n}{n+x^2-1} \Big)^{1/2}  \bigg ),
\]
see  Chung \citep{C46}.
So, an asymptotic bound on the tail probabilities for self-normalized sums implies an asymptotic bound on the tail probabilities for $T_n.$
Shao \cite{S97} established self-normalized large and moderate deviation principles  without any moment assumptions,
and
Shao \cite{S99} proved the following self-normalized Cram\'{e}r type moderate deviations: if  $\mathbb{E}|X_1|^{2+\rho}< \infty$ for some $ \rho \in (0, 1] ,$ then
\begin{equation}\label{Cramer01}
\frac{\mathbb{P}(W_n \geq x)}{1-\Phi \left( x\right)}=1+o(1)
\end{equation}
uniformly for  $0\leq x =o(n^{\rho/(4+2\rho)})$  as $n\rightarrow \infty.$ 
The 
later result
indicates that 
the normal range of $x$ for \eqref{Cramer01}   
on self-normalized sums can be much wider than that for classical moderate deviation expansion \eqref{Cramer00} on sums of i.i.d.\ r.v.'s.
 The expansion (\ref{Cramer01}) was further extended to independent but not necessarily identically distributed random variables by
Jing,   Shao and Wang \cite{JSW03}.
Their result 
implies the following precise asymptotic  normality under  finite $(2+\rho)$-th moments:
\begin{equation}\label{Cramer02}
\frac{\mathbb{P}(W_n \geq x)}{1-\Phi \left( x\right)}=\exp\Big\{   O\big( 1\big) \frac{ (1+x)^{2+\rho}}{n^{\rho/2}}   \Big\},
\end{equation}
uniformly for   $0\leq x =o(  \sqrt{n}) $ as $n \rightarrow \infty$.
Moderate deviation results of types  (\ref{Cramer01}) and (\ref{Cramer02})  play  an important role in statistical inference of means since
in practice one usually does not know the variance $\sigma^2.$
Even when the later can be estimated, it is still advisable to
 use self-normalized sums to obtain more precise results.
 Due to these significant advantages, the limit theory for self-normalized sums
 attracts more and more attention.
  Gin\'{e},  G\"otze and Mason \cite{GGM97} gave  a necessary and sufficient condition   for
the asymptotic normality of self-normalized partial sums.
Cs\"{o}rg\H{o}, Szyszkowicz and Wang \cite{CSW03} established Donsker's theorem. 
For various moderate and large deviations results for self-normalized sums, we refer to, for instance,
Jing,   Shao and Wang \cite{JSW03},
Liu, Shao and Wang \cite{LSW13},
de la Pe\~{n}a,   Lai and Shao \cite{DLS09},
Shao and Wang  \cite{SW13} and Shao \cite{S18}.
Dembo and Shao \cite{DS06} and Liu and Shao \cite{LS13} studied Hotelling's $T^2$-statistic.

The moderate deviation  theory for self-normalized sums of independent random variables
has been studied in depth. However, there are only a few results for dependent random variables.
Chen, Shao, Wu and Xu  \cite{CSWX16}   established  self-normalized Cram\'{e}r type moderate deviations for $\beta$-mixing sequences and functional  dependent sequences  (see  Wu \cite{W05} for the definition of functional dependent sequences).
Fan,  Grama, Liu and Shao \cite{FGLS17} gave two  self-normalized Cram\'{e}r type moderate deviation results
for martingales.  For a closely related topic, that is,  exponential inequalities for self-normalized martingales,
we refer to de la Pe\~{n}a \citep{D99} and Bercu and Touati \citep{BT08}.
The main purpose of this paper is to establish self-normalized Cram\'{e}r type moderate deviations for general stationary sequences.
We deduce also a self-normalized moderate deviation principle  and a Berry-Esseen bound.

The paper is organized as follows. Our main results are stated and discussed in Section \ref{sec2}.
The applications   are given in Section  \ref{sec3}.
Proofs of theorems are deferred to Section \ref{sec4}.

All over the paper, $c$ and $C$, possibly enabled with indices (arguments), denote constants depending only on the previously introduced
constants and on its indices (arguments).
Their values may change on every occurence.
For two positive real sequences $(a_n)_{i\geq1}$ and $(b_n)_{i\geq1},$
we write $a_n=O(b_n)$ if there exists a positive constant $C$
such that $ a_n  \leq C  b_n $ holds for all large $n$,
and $a_n=o(b_n)$ if $\lim_{n\rightarrow\infty}a_n/b_n =0.$
We also write $a_n \asymp b_n $ if $  a_n =O(b_n)$ and $b_n=O(a_n),$
and $a_n \sim b_n$ if $ \lim_{n\rightarrow \infty} a_n /b_n =1.$

\section{Main results}  \label{sec2}
\setcounter{equation}{0}

Assume that $(X_i)_{i \in \mathbb{Z}}$ is a stationary sequence of centered random variables, where $X_i=X_0 \circ T^i $ and
$T: \Omega \mapsto \Omega$ is a bijective bimeasurable transformation preserving the probability $\mathbb{P}$ on $(\Omega, \mathcal{F})$.
For a subfield $\mathcal{F}_0$ satisfying $\mathcal{F}_0 \subseteq T^{-1}(\mathcal{F}_0)$, let $\mathcal{F}_i= T^{-i}(\mathcal{F}_0).$
Assume  that $X_0$ is   $\mathcal{F}_0$-measurable, so that the sequence $(X_i)_{i \in \mathbb{Z}}$ is adapted to the filtration $(\mathcal{F}_i)_{i\in \mathbb{Z}}$.  

Denote by  $\lfloor a \rfloor$ the integer part of the real $a.$
Let $m \in  [1, n]$ and $ k=\lfloor n/m   \rfloor,$
where $m$ may depend on $n.$
Define
$$H_j= \{i: m(j-1)+1  \leq i \leq mj  \}, \ \ \ \ 1 \leq j \leq k.$$
Consider the block sums $S_{j}^\circ= \sum_{i \in H_j} X_i$, and the block self-normalized sums
$$ W^\circ_n= \frac{\sum_{j=1}^k S_{j}^\circ  }{ V_k^\circ   }, \ \    \ \ \ \ \ \ \textrm{where } \ \   (V_k^\circ)^2 =  \sum_{j=1}^k (S_{j}^\circ)^2   .$$
In particular, when $m=1,$ the  block self-normalized sum $W^\circ_n$ becomes self-normalized sum  $W_n.$
We also denote the $\mathbb{L}^\infty$-norm of $X$ by $\|X \|_\infty$, that is
$\|X \|_\infty = \inf\{u: \mathbb{P}(|X|> u) =0\}.$
For any $1\leq m\leq   n$, set
\begin{eqnarray}
&&  \varepsilon_m=    \frac {1}{n^{  1/2}m^{1/\rho} \sigma^{2/\rho +1  } }  \bigg \|    \mathbb{E}[ |S_m|^{2+\rho} | \mathcal{F}_{0}]      \bigg \|_\infty^{1/\rho},
\label{grmma1} \\
&&\gamma_m=  \frac{1}{m^{1/2} \sigma } \sum_{j=1}^{\infty} \frac{1}{ j^{3/2}} \Big \|\mathbb{E}[ S_{mj} | \mathcal{F}_0]\Big \|_\infty   \label{grmman}
\end{eqnarray}
 and
\begin{eqnarray}
\delta_m^2= \frac{1}{m \sigma^2 }  \bigg \|\mathbb{E}[ S_m |\mathcal{F}_0]\bigg\|_\infty^2   +   \bigg \| \frac1{m\sigma^2}  \mathbb{E}[ S_m^2| \mathcal{F}_0]- 1\bigg \|_\infty  ,
\end{eqnarray}
where $\rho$ and $ \sigma$ are two positive constants.
   We are interested in the case where
\begin{eqnarray}\label{bhg03}
\max\{\varepsilon_m,\, \gamma_m, \,  \delta_m,   m/n   \} \rightarrow 0 \ \ \ \textrm{as}  \ n\rightarrow \infty.
\end{eqnarray}
We remark that $\delta_m \rightarrow0$ implies 
that $\frac1m \sum_{i=1}^m \mathbb{E}S_m^2 \to \sigma^2$ as  $n\rightarrow \infty$.

\begin{remark}\label{re01}
Let us comment on   condition (\ref{bhg03}).
\begin{enumerate}
 \item If  $ \big \|    \mathbb{E}[ |X_1|^{2+\rho} | \mathcal{F}_{0}]      \big \|_\infty < \infty,$   then, by convexity,
 we have $$   \|  \mathbb{E}[ |\frac1m S_m|^{2+\rho} | \mathcal{F}_{0}] \|_\infty  \leq  \frac1m  \sum_{i=1}^m \|  \mathbb{E}[ |X_i|^{2+\rho} | \mathcal{F}_{0}] \|_\infty  \leq   \big \|    \mathbb{E}[ |X_1|^{2+\rho} | \mathcal{F}_{0}]      \big \|_\infty$$
 and thus $\varepsilon_m= O(  m^{1+1/\rho}/n^{1/2})$ as $n \rightarrow \infty.$  In particular, the claim holds provided that $X_1$ is bounded, that is
  $ \big \|  X_1    \big \|_\infty < \infty. $

 \item If  $   \|  X_1     \|_\infty < \infty$ and $\delta_m \rightarrow 0$,   then we have
 $$ \big \|    \mathbb{E}[ |S_m|^{2+\rho} | \mathcal{F}_{0}]      \big \|_\infty \leq m^{\rho}\big \|  X_1    \big \|_\infty^\rho \big \|    \mathbb{E}[  S_m ^{2 } | \mathcal{F}_{0}]      \big \|_\infty =O(m^{1+\rho}). $$
 Therefore, it holds  $\varepsilon_m= O(  m /n^{1/2})$ as $n \rightarrow \infty.$

 \item Assume   $ \big \|    \mathbb{E}[ |S_m|^{2+\rho} | \mathcal{F}_{0}]      \big \|_\infty=O(m^{1+\rho/2})$ as $ m\rightarrow \infty.$ Then it is easy to see that  $\varepsilon_m= O(  \sqrt{m/n})$ as $n \rightarrow \infty.$ In particular, if $(X_i, \mathcal{F}_i)_{i \in \mathbb{Z}}$ is a  martingale difference sequence satisfying $\big \| \mathbb{E}[ |X_1|^{2+\rho} | \mathcal{F}_{0}]  \big \|_\infty < \infty$,  then, by Theorem 2.1 of Rio   \cite{R09}, it is easy to see that
     $$    (  \mathbb{E}[ |S_m|^{2+\rho} | \mathcal{F}_{0}] )^{2/(2+\rho)}   \leq (1+\rho) \sum_{k=1}^{m}  (  \mathbb{E}[ |X_i|^{2+\rho} | \mathcal{F}_{0}] )^{2/(2+\rho)}  \leq (1+\rho) \big \| \mathbb{E}[ |X_1|^{2+\rho} | \mathcal{F}_{0}]  \big \|_\infty^{2/(2+\rho)}m  $$
     a.s., which leads to
     $$ \big \|    \mathbb{E}[ |S_m|^{2+\rho} | \mathcal{F}_{0}]      \big \|_\infty=O(m^{1+\rho/2}) \ \ \  \emph{and} \ \ \ \varepsilon_m= O(  \sqrt{m/n}) $$
     as $n\to\infty$, and
     $$ \gamma_m= 0 \ \ \  \emph{and} \ \ \  \delta_m^2=    \Big \| \frac1{m\sigma^2}  \sum_{i=1}^m \mathbb{E}[ X_i^2| \mathcal{F}_0]- 1\Big \|_\infty.$$

 \item  Dedecker et al.\ \cite{DMPU09} introduced the following two conditions for stationary sequences:
\begin{description}
  \item[(A1)] The following sum is finite:
  \begin{equation}\label{cond01}
\sum_{n=1}^{\infty} \frac{1}{n^{3/2}}   \Big \|\mathbb{E}[ S_n| \mathcal{F}_0] \Big \|_\infty < \infty.
\end{equation}
  \item[(A2)] There exists a positive constant $\sigma  $ such that
\begin{equation}\label{cond02}
 \lim_{n\rightarrow\infty} \Big \| \frac 1 n  \mathbb{E}[ S_n^2| \mathcal{F}_0]- \sigma^2\Big \|_\infty=0.
\end{equation}
\end{description}
Clearly,  under conditions (A1) and (A2),  by Lemma 29 of Dedecker et al.\ \cite{DMPU09},  it holds that $\max\{ \gamma_m, \,   \delta_m  \} \rightarrow 0$  for any sequence $ m = m(n)$   such that $m\rightarrow \infty$ and $m /n  \rightarrow 0$  as $n\rightarrow \infty$.

\end{enumerate}
\end{remark}

For any sequence of small positive numbers $(\varepsilon_m)_{m\geq 1},$  let  $\widehat{\varepsilon}_m(x, \rho)$  be a function of $\varepsilon_m, x$ and $\rho$ defined as follows
 \begin{equation}\label{defepsilonm}
 \widehat{\varepsilon}_m(x, \rho) = \frac{ \varepsilon_m^{ \rho(2-\rho)/4 } }{1+ x^{\rho(2+\rho)/4}} .
 \end{equation}
The following  theorem  gives a self-normalized Cram\'{e}r type moderate deviation result  for  stationary sequences.
\begin{theorem}\label{th1}
 Assume that there exists $\rho \in (0,    1 ]$  such that $\max\{\varepsilon_m,\, \gamma_m, \,   \delta_m, m/n  \} \rightarrow 0$ as $n\rightarrow \infty.$
\begin{description}
  \item[\textbf{[i]}] If $\rho \in (0, 1)$, then   there exists an absolute constant $\alpha_{\rho   }   >0 $ such that for all $0\leq x \leq \alpha_{\rho   } \min\{\varepsilon_m^{-1},\, \sqrt{n/m} \}  , $
\begin{eqnarray*}
\Bigg| \ln \frac{\mathbb{P}( W_n^\circ  \geq x   )}{1-\Phi \left(  x\right)}  \Bigg|  &\leq&  C_{\rho} \Bigg( x^{2+\rho}  \varepsilon_m^\rho+ x^2 \Big(\delta_m^2   +\gamma_m |\ln \gamma_m|+\frac{m}n \Big)   \\
 &&\ \ \ \ \ \   +\ (1+x)\Big(  \delta_m +\gamma_m |\ln \gamma_m| +\varepsilon_m ^\rho + \widehat{\varepsilon}_m(x, \rho)+\sqrt{\frac{m}n} \Big) \Bigg)  ,
\end{eqnarray*}
where  $C_{ \rho }$ depends only on $ \rho.$

 \item[\textbf{[ii]}] If $\rho =1$, then   there exists an absolute constant $\alpha  >0 $ such that for all $0\leq x \leq \alpha  \min\{\varepsilon_m^{-1},\, \sqrt{n/m} \}  , $
\begin{eqnarray*}
\Bigg| \ln \frac{\mathbb{P}( W_n^\circ  \geq x   )}{1-\Phi \left(  x\right)}  \Bigg|  &\leq&  C  \Bigg( x^{3}  \varepsilon_m + x^2 \Big(\delta_m^2   +\gamma_m |\ln \gamma_m|+\frac{m}n \Big)   \\
 &&\ \ \ \ \    +\ (1+x)\Big(  \delta_m +\gamma_m |\ln \gamma_m| + \varepsilon_m |\ln \varepsilon_m| + \widehat{\varepsilon}_m(x, 1) + \sqrt{\frac{m}n} \Big) \Bigg) .
\end{eqnarray*}
\end{description}
In particular, the last two inequalities   imply that, for any $\rho \in (0,1],$
\begin{eqnarray}\label{thls}
 \frac{\mathbb{P}( W_n^\circ  \geq x    )}{1-\Phi \left(  x\right)}   = 1 +o(1)
\end{eqnarray}
uniformly  for   $\displaystyle  0 \leq x = o \big(  \min \big\{ \varepsilon_m^{-\rho/(2+\rho)}, \, \delta_m^{-1}  \, , (\gamma_m |\ln \gamma_m|)^{-1/2}, \sqrt{n/m} \big\}  \big)$ as $n \rightarrow \infty.$
Moreover, the same results hold  with $\displaystyle \frac{\mathbb{P}( W_n^\circ  \geq x  )}{1-\Phi \left(  x\right)}$   replacing by $\displaystyle \frac{\mathbb{P}(W_n^\circ \leq -x )}{ \Phi \left(-  x\right)}$.
\end{theorem}

\begin{remark}\label{remark2.1}
Let us comment on the results of Theorem  \ref{th1}.
\begin{enumerate}
\item  The range of validity of (\ref{thls}) can be very large. For instance, if   $\big \|\mathbb{E}[ |S_n|^{2+\rho} |\mathcal{F}_0]\big \|_\infty =O(n^{1+\rho/2 }),$  $\big \|\mathbb{E}[ S_{n} | \mathcal{F}_0]\big \|_\infty=O(1),$ and $\big\| \frac1 n  \mathbb{E}[ S_n^2| \mathcal{F}_0]- \sigma^2\big\|_\infty=O\big( \frac 1 n\big) $
as $n\rightarrow \infty,$ then
 $$\varepsilon_m= O(  \sqrt{m/n}),\ \ \ \ \   \gamma_m, \delta_m= O(  \sqrt{1/m  }) .$$
 With $m=\lfloor n^{2\rho /(2+3\rho)} \rfloor,$  equality  (\ref{thls}) holds uniformly  for $ 0 \leq x = o(n^{\rho/(4+6\rho)} /\sqrt{\ln n} )$ as $n \rightarrow \infty.$ The last range coincides with   the classical range, up to a term   $\sqrt{\ln n}$, when applied for block self-normalized sums of i.i.d.\,random variables, that is $ 0 \leq x = o(k^{\rho/(4+2\rho)}   ).$ See Remark 1 of Shao \cite{S99}.

\item If $(X_i, \mathcal{F}_i)_{i \in \mathbb{Z}}$ is a  martingale difference sequence satisfying $\|\mathbb{E}[ |X_1|^{2+\rho} | \mathcal{F}_{0}] \|_\infty  < \infty$,  then Theorem \ref{th1} gives a block self-normalized Cram\'{e}r type moderate deviation result,
    with
  $$\varepsilon_m = O \big( \sqrt{ m/n } \big),\ \   \displaystyle \gamma_m=0 \ \   \textrm{and} \ \   \delta_m^2=  \Big \| \frac1 {m \sigma^2} \sum_{i=1}^m \mathbb{E}[ X_i^2| \mathcal{F}_0]- 1\Big \|_\infty $$
    as $n\rightarrow \infty, $ which extends the main result of Fan et al.\,\cite{FGLS17} to block self-normalized  martingales.
    Furthermore, if $ \|  \mathbb{E}[ X_i^2| \mathcal{F}_0] -\sigma^2\|_\infty  \leq Ci^{-\theta}$  for some positive  constants $C$ and $\theta,$
    then  we have
  \begin{displaymath}
 \delta_m ^2=  \left\{ \begin{array}{ll}
  O( m^{-1} ), & \textrm{\ \ \  if $\theta > 1$,}\\
 O( m^{-1}  \ln m  ),  & \textrm{\ \ \ if $\theta = 1$,} \\
 O( m^{-\theta} ), & \textrm{\ \ \  if $\theta \in (0, 1)$.}
\end{array} \right.
\end{displaymath}
Taking
\begin{displaymath}
 m=  \left\{ \begin{array}{ll}
   \lfloor n^{\rho/(2+2\rho)} \rfloor, & \textrm{\ \ \  if $\theta \geq 1$,}\\
 \lfloor n^{ \rho /(  \rho+\theta(2+\rho))} \rfloor, & \textrm{\ \ \  if $\theta \in (0, 1)$,}
\end{array} \right.
\end{displaymath}
we have the following results:
\begin{description}
\item[\textbf{[i]}] If $\rho \in (0, 1)$, then (\ref{thls}) holds for  $0\leq x =o(n^{ \theta \rho /(2\rho+2\theta(2+\rho))})$.
\item[\textbf{[ii]}] If $\rho =1$, then (\ref{thls}) holds for  $0\leq x =o(n^{ \rho /(4+4\rho)}/ \ln n).$
\item[\textbf{[iii]}] If $\rho >1$, then (\ref{thls}) holds for  $0\leq x =o(n^{ \rho /(4+4\rho)})$.
\end{description}

\item  Besides block self-normalized sums, we can also consider the interlacing self-normalized sums. Let  $\alpha \in (0, 1)$
and
$m=\lfloor n^\alpha \rfloor, \ k =\lfloor n/(2m)  \rfloor $ (instead of $\lfloor n/m \rfloor $ considered before) and
$$B_j= \Big\{i: 2m(j-1)+1  \leq i \leq 2mj-m  \Big \}, \ \ \ \ 1 \leq j \leq k.$$
Let $Y_j^*= \sum_{l \in B_j} X_l, \ (V_k^*)^2= \sum_{j=1}^k (Y_j^*)^2$ and write
$$I_n^*= \frac{ \sum_{j=1}^k Y_j^* }{  V_k^* }$$
for the interlacing self-normalized sum.
Clearly, Theorem \ref{th1} also holds for interlacing self-normalized sums $I_n^*$,
 with $  \mathbb{E}[ \, \cdot\, |\mathcal{F}_0]$  and $W_n^\circ$ replaced respectively by $ \mathbb{E}[ \, \cdot\, |\mathcal{F}_{-m}]$ and
$I_n^*$.
Such type of results for $\beta$-mixing and some functional  dependent sequences have been considered  by  Chen et al.  \cite{CSWX16}.
\end{enumerate}
\end{remark}

The following self-normalized moderate deviation principle (MDP) result  is a  consequence of Theorem \ref{th1}.
\begin{corollary}\label{co0}
Assume the condition of Theorem \ref{th1}.   Let $(a_n)_{n\geq1}$ be any sequence of real numbers satisfying $a_n \rightarrow 0$ and $a_n  \min\{\varepsilon_m^{-1},\, \sqrt{n/m} \}\rightarrow \infty$
as $n\rightarrow \infty$.  Then,  for each Borel set $B \subset \mathbb{R}$,
\begin{eqnarray}
- \inf_{x \in B^o}\frac{x^2}{2   }  \leq   \liminf_{n\rightarrow \infty}  a_n^2 \ln \mathbb{P}\Big(a_n W_n^o      \in B \Big)
  \leq  \limsup_{n\rightarrow \infty} a_n^2\ln \mathbb{P}\Big(a_n  W_n^o     \in B \Big) \leq  - \inf_{x \in \overline{B}}\frac{x^2}{2  }\,,   \label{MDP}
\end{eqnarray}
where $B^o$ and $\overline{B}$ denote the interior and the closure of $B$, respectively.
\end{corollary}

In the i.i.d.\,case, $W_n^o$ is a self-normalized sum of $k$ i.i.d.\ random variables. 
According to the classical result of Jing, Shao and Wang \cite{JSW03}, the MDP holds for $0\leq x =o(k^{1/2}).$
Since $k=\lfloor n/m\rfloor$, the last range reads also as $0\leq x  =o( \sqrt{n/m} )$. Notice that $ \varepsilon_m^{-1}$ is of order $\sqrt{n/m}$.
Thus, the convergence rate of $a_n$ in the  last corollary cannot be improved even for  i.i.d.\ random variables.

Theorem \ref{th1} also implies the following self-normalized Berry-Esseen bound  for  stationary sequences.
\begin{corollary}\label{co01}
Assume the condition of Theorem \ref{th1}.  Then, for $\rho \in (0,    1 ]$,
\begin{eqnarray*}
\sup_{ x   } \Big|\mathbb{P}( W_n^\circ \leq x   )  - \Phi \left(  x\right) \Big|   \leq   C_{\rho}  \Big(  \delta_m +\gamma_m |\ln \gamma_m| +\varepsilon_m^{ \rho(2-\rho)/4 }+\sqrt{\frac{m}n} \Big)  ,
\end{eqnarray*}
where  $C_{ \rho }$ depends only on $ \rho.$
\end{corollary}


\section{Applications} \label{sec3}
\setcounter{equation}{0}
 In this section, we present some applications of our results, including
 $\phi$-mixing type sequences,
 contracting Markov chains,  expanding maps
  and
  confidence intervals.

\subsection{$\phi$-mixing type sequences}
Let $Y$ be a random variable with values in a Polish space $\mathcal{Y}.$
If $\mathcal{M}$ is a $\sigma$-field,
the $\phi$-mixing coefficient between $\mathcal{M}$ and $\sigma(Y)$ is defined by
\begin{equation}\label{phi}
\phi(\mathcal{M}, \sigma(Y))=\sup_{A \in \mathfrak{B}(\mathcal{Y})} \Big\| \mathbb{P}_{Y|\mathcal{M}}(A) -\mathbb{P}_Y(A) \Big\|_\infty.
\end{equation}
For  a sequence of  random variables $(X_i)_{i \in \mathbb{Z}}$ and a positive integer $m,$ denote
\[
\phi_m(n)=\sup_{i_m> ...> i_1\geq n}\phi(\mathcal{F}_0, \sigma(X_{i_1},...,X_{i_m} )),
 \]
 and let $\phi(k)=\lim_{m\rightarrow \infty}\phi_m(k)$ be the usual $\phi$-mixing coefficient.
Under  the following condition
\begin{equation}\label{phic01}
 \sum_{k\geq 1}k^{-1/2}\phi_1(k) < \infty  \ \   \ \ \textrm{and }\ \ \ \ \lim_{k \rightarrow \infty} \phi_2(k)=0,
\end{equation}
 Dedecker et al.\ \cite{DMPU09} obtained a MDP for standardized sums of bounded $\phi$-mixing random variables.
See also Gao \cite{G96} for an earlier version of MDP under the condition $ \sum_{k\geq1} \phi(k) < \infty$
which is stronger than (\ref{phic01}).
Denote
 \begin{eqnarray*}
 \eta_{1, n}&=& \sup_{k\geq n} \| \mathbb{E}[X_k |\mathcal{F}_0]\|_\infty, \\
  \eta_{2, n} &=& \sup_{k, l\geq n} \| \mathbb{E}[X_kX_l |\mathcal{F}_0] -\mathbb{E}[X_kX_l ] \|_\infty.
  \end{eqnarray*}
Clearly, when the random variable $ X_0 $ is   bounded, it holds that $\eta_{1, n}=O(\phi_1(n))$ and $\eta_{2, n}=O(\phi_2(n))$ as $n\rightarrow \infty$.

From Theorem \ref{th1} we obtain the following  self-normalized   Cram\'{e}r type moderate deviation expansion
with depending structure defined by $\eta_{1, n}$ and $\eta_{2, n}.$

\begin{proposition}\label{pro3.3}   Assume that $  \|X_0\|_\infty   <\infty,$    $$  \sigma^2:= \sum_{k=-\infty}^{\infty}\mathbb{E}[X_0X_k ] >0 \ \ \ \textrm{and}  \ \ \ \max_{i=1,2 }\{   \eta_{i, n}  \} =O(  n^{-\beta}),\ \ n\rightarrow \infty ,$$ for some constant $\beta>1.$
\begin{description}
  \item[\textbf{[i]}] If $\beta\geq 3/2$,
then     (\ref{thls}) with $m=\lfloor n^{2/7} \rfloor$  holds uniformly  for $ 0 \leq x = o(n^{1/14} /\sqrt{\ln n} )$ as $n \rightarrow \infty.$

 \item[\textbf{[ii]}] If   $\beta \in (1, 3/2),$
then   (\ref{thls}) with $m=\lfloor n^{  1/(3\beta-1)} \rfloor$   holds uniformly  for $0 \leq x = o(n^{(\beta-1)/(6\beta-2)})$ as $n \rightarrow \infty.$

\item[\textbf{[iii]}]
Assume $m :=m(n)\rightarrow \infty$  and $n^{ 1/2} /m \rightarrow \infty$  as $n \rightarrow \infty.$  Let $(a_n)_{n\geq1}$ be any sequence   of real numbers such that $a_n \rightarrow 0$ and   $a_n   n^{ 1/2} /m  \rightarrow \infty$
as $n\rightarrow \infty$. Then (\ref{MDP}) holds.
\end{description}
\end{proposition}

By point 3 of Remark \ref{re01},   if    $\mathbb{E} |S_n|^{2+\rho}    =O(n^{1+\rho/2}) $  for some $\rho>0,$  then point [iii] of Proposition \ref{pro3.3} can be
further improved. Indeed, in this case,   (\ref{MDP}) holds
for any $m \rightarrow \infty,$ and any sequence of real numbers $(a_n)_{n\geq1}$ such that $a_n \rightarrow 0$ and $a_n   \sqrt{n/m}  \rightarrow \infty$
as $n\rightarrow \infty$.


\subsection{Contracting Markov chains}\label{ap3.3}
Let $(Y_n)_{n\geq 0} $ be a stationary Markov chain of  bounded random variables  with invariant measure $\mu$ and transition kernel $K.$
Denote by $\| \cdot \|_{\infty, \mu} $ the essential norm with respect to $\mu.$ Let $\Lambda_1$ be the set of $1$-Lipschitz
functions. Assume that the Markov chain satisfies the  following condition:
\begin{description}
  \item[(B)] \, There exist two constants $C>0$ and $\rho \in (0, 1)$ such that \[\sup_{g \in \Lambda_1} \|  K^n(g) - \mu(g)  \|_{\infty, \mu} \leq C \rho^n   \]
    and  for any $ g , g' \in \Lambda_1$ and any $m\geq 0,$
  \[  \lim_{n \rightarrow \infty}   \Big\|  K^n\big(g'K^m(g) \big) - \mu\big(g'K^m(g)\big)  \Big\|_{\infty, \mu}   =0. \]
\end{description}

Denote by   $\mathcal{L}$   the class of functions
$f:  \mathbf{R} \mapsto \mathbf{R}$ such that
\begin{equation}
\label{fsddfs}
|f(x) -f(y)| \leq h(|x-y|),
\end{equation}
where $h$ is a  concave and non-decreasing function
satisfying
\begin{equation}
\label{fsddsdfs}
 \int_0^1 \frac{h(t)}{ t \sqrt{|\ln t| }}  dt < \infty,
\end{equation}
see \cite{DMPU09}.
Clearly, inequality  (\ref{fsddsdfs}) holds if $h(t) \leq c |\ln(t)|^{-\gamma}$ for some constants $c>0$ and  $\gamma > 1/2.$
In particular,  $\mathcal{L}$ contains the class of $\alpha$-H\"{o}lder continuous  functions from $[0, 1]$ to $\mathbf{R}$, where $\alpha \in (0, 1].$

Dedecker et al.\ \cite{DMPU09} proved a   MDP 
for the sequence
 \begin{equation} \label{fds}
 X_n=f(Y_n) - \mu(f)
 \end{equation}
under the condition that the function $f$ belongs to the class   $\mathcal{L}.$
The following proposition gives an extension of the  MDP to \textit{self-normalized} sums
$ W^\circ_n= \frac{\sum_{j=1}^k S_{j}^\circ  }{ V_k^\circ   }, $
where $S_{j}^\circ= \sum_{i \in H_j} X_i$ and $ (V_k^\circ)^2 =  \sum_{j=1}^k (S_{j}^\circ)^2   .$
\begin{proposition}\label{ap02s}
Assume that the stationary Markov chain $(Y_n)_{n\geq 0}$ satisfies condition (B), and let $X_n$
be defined by (\ref{fds}), with $f$ belonging to $\mathcal{L}$. Assume $m :=m(n)\rightarrow \infty$  and $n^{ 1/2} /m \rightarrow \infty$  as $n \rightarrow \infty.$ Let $a_n$ be any sequence   of real numbers such that $a_n \rightarrow 0$ and $a_n   n^{ 1/2} /m  \rightarrow \infty$
as $n\rightarrow \infty$. If
$$\sigma^2:=\sigma^2(f)= \mu\big ((f-\mu(f))^2 \big) + 2 \sum_{n\geq1} \mu \big(   K^n(f) \; (f-\mu(f)) \big) >0,  $$
then (\ref{MDP}) holds.
\end{proposition}
 \emph{Proof.}
 By Lemma 15 of Dedecker et al.\ \cite{DMPU09}, it is easy to see that $X_1$ is bounded:  
$\| X_1\|_{\infty, \mu} \leq h(C \rho) $ with $h$ defined by (\ref{fsddfs}).
 Then by point 2 of Remark \ref{re01},
we have  $\varepsilon_m= O(  m /n^{1/2})$ as $n \rightarrow \infty.$
The conditions of Proposition \ref{ap02s} imply the conditions (\ref{cond01}) and (\ref{cond02}): see
the proof of Proposition 14 in Dedecker et al.\ \cite{DMPU09}.
Hence, by point 4 of Remark \ref{re01},
 the conditions of  Proposition \ref{ap02s}
 imply the conditions of Corollary  \ref{co0}, thus Proposition \ref{ap02s} follows. 
 \hfill\qed

Furthermore, assume that the Markov chain satisfies the following condition  which is stronger than condition (B).
\begin{description}
  \item[(C)] \, There exist two constants $C>0$ and $\rho \in (0, 1)$ such that \[\sup_{g \in \Lambda_1} \|  K^n(g) - \mu(g)  \|_{\infty, \mu} \leq C \rho^n   \]
    and   for any $m\geq 0,$
  \[    \sup_{g , g' \in \Lambda_1} \Big\|  K^n\big(g'K^m(g) \big) - \mu\big(g'K^m(g)\big)  \Big\|_{\infty, \mu}  \leq C \rho^n. \]
\end{description}
  Then we have  the following  self-normalized   Cram\'{e}r type moderate deviation expansion.
\begin{proposition} \label{pr3.3s}
Assume that the stationary Markov chain $(Y_n)_{n\geq 0}$ satisfies condition  (C), and let $X_n$
be defined by (\ref{fds}).
Assume $f \in \mathcal{L},$
\[
\sigma^2:=\sigma^2(f)=\mu \Big((f-\mu(f))^2 \Big) + 2 \sum_{n>0} \mu \Big(K^n(f)\cdot (f-\mu(f))\Big)  >0
\]
and, for  some constant $\beta>1,$
 \begin{equation} \label{dfs456}
 h(   \rho^n ) = O(  n^{-\beta}) , \ \ \ \ n\rightarrow\infty,
 \end{equation}
 where  $h$ is defined by (\ref{fsddfs}).
\begin{description}
  \item[\textbf{[i]}]  If $\beta \geq 3/2$, then  (\ref{thls})    with $m=\lfloor n^{2/7} \rfloor$ holds uniformly  for $ 0 \leq x = o(n^{1/14} /\sqrt{\ln n} )$ as $n \rightarrow \infty.$
  \item[\textbf{[ii]}]   If $\beta \in (1, 3/2)$, then  (\ref{thls}) with $m=\lfloor n^{  1/(3\beta-1)} \rfloor$  holds uniformly  for $ 0 \leq x = o(n^{(\beta-1)/(6\beta-2)}   )$ as $n \rightarrow \infty.$
\end{description}
\end{proposition}

Notice that if $g(t) \leq c |\ln(t)|^{-\beta}$ for some constants $c>0$ and  $\beta > 1,$ then   (\ref{dfs456}) is satisfied. \\
\emph{Proof.} From the proof of Propositions 14  of  \cite{DMPU09}, it is easy to see that
 $$ \max_{i=1,2}\{\eta_{i, n}  \} =O \big(  h(C \rho^n ) \big)  ,$$
 where $C$ is given by condition (C) and $h$ is defined by (\ref{fsddfs}). Notice that $C\rho^n \leq \rho^{n/2}$ for  $n $ large enough.  Hence, Proposition  \ref{pr3.3s} is a simple consequence of Proposition \ref{pro3.3}.  \hfill\qed

\subsection{Expanding maps}
Dedecker et al.\ \cite{DMPU09} have obtained the MDP for expanding maps. Here we show  that our results can also be applied to
 expanding maps for getting self-normalized MDP and Cram\'{e}r type moderate deviations.

Let $T$ be a map from $[0, 1]$ to $[0, 1]$ preserving a probability $\mu$ on $[0, 1],$ and denote
$$X_k=f\circ T^{n-k+1}-\mu(f),$$
for any function $f \in L^2([0, 1], \mu)$.
Let
$ W^\circ_n= \frac{\sum_{j=1}^k S_{j}^\circ  }{ V_k^\circ   }, $
where $S_{j}^\circ= \sum_{i \in H_j} X_i$ and $ (V_k^\circ)^2 =  \sum_{j=1}^k (S_{j}^\circ)^2   .$
Denote by $\mathcal{BV}$  the class of bounded variation functions from $[0, 1]$ to $\mathbf{R}$. For any $f \in \mathcal{BV},$ denote by
$\|df\|$ the total variation norm of the measure $df:  \|df\|=\sup\{ \int g df, \|g\|_\infty \leq 1 \}.$  A Markov kernel $K$  is said to be
$\mathcal{BV}$-contracting if there exist two constants $k >0$ and $\rho \in [0, 1)$ such that
\begin{equation}\label{thkds02q}
\|d K^n(f) \| \leq k \rho^n \|df\|.
\end{equation}
Define the Perron-Frobenius operator $K$ from $L^2([0, 1], \mu)$ to $L^2([0, 1], \mu)$ via the equality
\begin{equation}\label{thkineq}
   \int_0^1 (K h)(x)f(x) \mu(dx) =\int_0^1 h(x)  ( f\circ T)(x) \mu(dx).
\end{equation}
The map $T$ is said to be $\mathcal{BV}$-contracting if its Perron-Frobenius operator is $\mathcal{BV}$-contracting.
We have the following corollary for the self-normalized sum $W_n^\circ.$

\begin{proposition}\label{pssfr4}
Assume that $T$ is $\mathcal{BV}$-contracting, $f \in \mathcal{BV}$ and
$\sigma^2:= \mu\big((f-\mu(f))^2\big) + 2 \sum_{n\geq2} \mu\big( f\circ T^n\cdot (f-\mu(f)) \big) >0  $.
\begin{description}
 \item[\textbf{[i]}] Let $m=\lfloor n^{2/7} \rfloor$. Equality (\ref{thls}) holds uniformly  for $ 0 \leq x = o(n^{1/14} /\sqrt{\ln n} )$ as $n \rightarrow \infty.$
 \item[\textbf{[ii]}] Assume $m :=m(n)\rightarrow \infty$  and $n^{ 1/2} /m \rightarrow \infty$  as $n \rightarrow \infty.$ Let $(a_n)$ be any sequence of real numbers such that $a_n \rightarrow 0$ and $a_n   n^{ 1/2} /m  \rightarrow \infty$
as $n\rightarrow \infty$.  Then (\ref{MDP}) holds.
\end{description}
\end{proposition}
 \emph{Proof.}
 Let $(Y_i)_{i\geq 1}$ be the Markov chain with transition kernel $K$ and invariant measure $\mu$
 in the stationary regime.
 Using equality (\ref{thkineq}), it is easy to see that $(Y_0, ..., Y_n)$ is distributed as $(T^{n+1},...,T).$
 Assume that $f \in \mathcal{BV}$. Since $K$ is $\mathcal{BV}$-contracting, by  the proof of Corollary 18 of \cite{DMPU09}, we have
 $$\|\mathbb{E}[ X_k | Y_0] \|_\infty  
 \leq C\rho^k \|df \|$$
 and, for any $l>k \geq 0,$
 $$\|\mathbb{E}[ X_kX_l | Y_0] - \mathbb{E}[ X_kX_l  ]  \|_\infty \leq C(1+C) \rho^k \|df \|^2.$$
By an argument similar to the proof of  Proposition \ref{pro3.3},  Proposition  \ref{pssfr4}  follows.
 \hfill\qed

\subsection{Application to confidence intervals}
Consider the problem of constructing  confidence intervals for the mean value $\mu$ of the stationary sequence $(\zeta_i)_{i\geq 1}$.
Let $X_i=\zeta_i-\mu, i\geq 1.$
Assume that $(X_i)_{i\geq 1}$ satisfies  the  conditions  (\ref{grmma1})-(\ref{bhg03}).
Let
$$T_n= \frac{\sum_{j=1}^k ( Y_j - m\mu ) }{\sqrt{\sum_{j=1}^k (Y_j - \overline{Y}_j)^2}}, $$
where  $Y_j= \sum_{i=1}^m \zeta_{ m(j-1)+i}, \  1\leq j \leq k,$ and $\overline{Y}_j=k^{-1}\sum_{j=1}^k Y_j.$
\begin{proposition}\label{c0kls} Let $\kappa_n \in (0, 1).$    Assume that  $\kappa_n \rightarrow  0$  and
\begin{eqnarray}\label{keldet}
 \big| \ln \kappa_n \big| =o \Big(   \min \big\{ \varepsilon_m^{-2 }, \, n/m  \big\}  \Big),\ \ \  n\rightarrow \infty.
\end{eqnarray}
Let $\Delta_n=\frac{ \sqrt{2  | \ln (\kappa_n/2)  |}  }{km } \sqrt{ \sum_{j=1}^k (Y_j - \overline{Y}_j)^2}$.
Then $[A_n,B_n]$ with
\begin{eqnarray*}
A_n=\frac{\sum_{j=1}^k Y_j }{k m}  -\Delta_n, \quad 
B_n=\frac{\sum_{j=1}^k Y_j }{k m}  +\Delta_n, 
\end{eqnarray*}
is a $1-\kappa_n$ confidence interval for $\mu$, for $n$ large enough.
\end{proposition}
\emph{Proof.}
 It is well known that for all $x\geq0,$
\[
\mathbb{P}\Big( T_n  \geq x \Big) = \mathbb{P}\Bigg(  \frac{\sum_{j=1}^k ( Y_j - m\mu  ) }{\sqrt{\sum_{j=1}^k (Y_j -  m\mu )^2}}  \geq x \Big(  \frac{k}{k-1}  \Big)^{1/2} \Big(\frac{k}{k+x^2-1} \Big)^{1/2}  \Bigg),
\]
see Chung \cite{C46}. The last  equality  and Theorem \ref{th1}  together implies that
\begin{equation} \label{tphisns4}
\frac{\mathbb{P}(T_n  \geq x)}{1-\Phi \left( x\right)}=\exp\Big\{ o(1)(1+ x)^2  \Big\}
\end{equation}
uniformly for $0\leq x=o\big(  \min \big\{   \varepsilon_m^{-1}, \,\sqrt{n/m} \big\} \big).$
Let $F(x)=1-(1-\Phi \left( x\right)  )\exp\{  o(1)(1+ x)^2 \}$.
  Notice that $$1-\Phi \left( x_n\right) \rightarrow \frac{1}{x_n\sqrt{2\pi}}e^{-x_n^2/2}= \exp\bigg\{-\frac{x_n^2}{2}\Big(1+\frac{2}{x_n^2}\ln (x_n\sqrt{2\pi})   \Big) \bigg\} ,\ x_n \rightarrow \infty. $$ Thus the upper $(\kappa_n/2)$-th quantile of the distribution function $F$ satisfies
$$ F^{-1}( \kappa_n/2) \rightarrow  \sqrt{2   | \ln (\kappa_n/2)  |} ,\ \ \ n \rightarrow \infty,$$ which, by (\ref{keldet}), is of order $o\big(  \min \big\{   \varepsilon_m^{-1}, \,\sqrt{n/m} \big\} \big).$
Then applying (\ref{tphisns4}) to $T_n$, we complete the proof of  Proposition \ref{c0kls}. \hfill\qed

By  (\ref{keldet}), a good choice of
the size $m$ is such that $R_n :=\min \big\{ \varepsilon_m^{-2 }, \, n/m  \big\}$  is  large enough,  so that $\kappa_n$ can be small enough.
A suitable choice is $m=\lfloor \ln n \rfloor;$ then, by Remark \ref{re01}, we have
\begin{displaymath}\displaystyle
 R_n=  \left\{ \begin{array}{ll}
   \frac{n}{\lfloor \ln n \rfloor^2} , & \textrm{\ \ \  if $\|  X_1     \|_\infty < \infty$,}  \\
   \frac{n}{\lfloor \ln n \rfloor^{2+2/\rho}}  ,  & \textrm{\ \ \ if $ \big \|    \mathbb{E}[ |X_1|^{2+\rho} | \mathcal{F}_{0}]      \big \|_\infty < \infty$,}  \\
 \frac{n}{\lfloor \ln n \rfloor}, & \textrm{\ \ \  if $\big \|    \mathbb{E}[ |S_m|^{2+\rho} | \mathcal{F}_{0}]      \big \|_\infty=O(m^{1+\rho/2})$.}
\end{array} \right.
\end{displaymath}

  Proposition \ref{c0kls}  
  uses a condition  on the $\mathbb{L}^\infty$-norm.  We should mention that Hannan's central limit theorem  (cf.\ Hannan \cite{H73})  holds under the condition on the $\mathbb{L}^{{2}}$-norm.
Accordingly, 
 a 
 confidential interval
 for linear regression  can be obtained via Hannan's  theorem   (cf.\ Caron and Dede \cite{CD18}), but with larger risk probability; the risk probability can be significantly improved,   
  using Cram\'{e}r type moderate deviations of Wu and Zhao \cite{WZ08} and  Cuny and Merlev\`{e}de \cite{CM14} on stationary sequences.
 Notice that the results of  \cite{WZ08} and   \cite{CM14} also hold when
 $X_i$ has finite $p$-th moments with $p > 2$.
 See also  Chen et al.\ \cite{CSWX16} for self-normalized Cram\'{e}r type moderate deviations for $\beta$-mixing sequences and functional  dependent sequences.

\section{Proofs of Theorems}\label{sec4}
\setcounter{equation}{0}

The proofs of our results are mainly based on the following lemmas which give some exponential deviation inequalities
for the partial sums of dependent random variables.

\subsection{Preliminary lemmas}\label{sec4.0}
Assume
on the probability space $(\Omega ,\mathcal{F},\mathbb{P})$
we are given a sequence of martingale differences $(\xi_i,\mathcal{F}_i)_{i=0,...,n}$,
where $\xi_0=0 $,  $\{\emptyset, \Omega\}=\mathcal{F}_0\subseteq ...\subseteq \mathcal{F}_n\subseteq
\mathcal{F}$ are increasing $\sigma$-fields. Define
\begin{equation}
M_{0}=0,\ \ \ \ \ M_k=\sum_{i=1}^k \xi_i,\quad k=1,...,n.  \label{xk}
\end{equation}
Let $[M]_n$  and  $\left\langle M\right\rangle_n $  be respectively the squared variance and the conditional variance   of the
martingale $M=(M_k,\mathcal{F}_k)_{k=0,...,n}$, that is
\begin{equation}\label{quad}
 [ M ]_0=0,\ \ \   [ M ]_k=\sum_{i=1}^k \xi_i^2, \ \ \
\left\langle M\right\rangle _0=0,\ \ \   \left\langle M\right\rangle _k=\sum_{i=1}^k\mathbb{E}[\xi_i^2|\mathcal{F}
_{i-1}],\quad k=1,...,n.
\end{equation}
Assume   the following two conditions:
\begin{description}
\item[(C1)]  There exist   $\epsilon_n \in (0, \frac12]$ and $\rho \in (0, +\infty)$ such that
\[
\mathbb{E}[|\xi_{i}| ^{2+\rho}  | \mathcal{F}_{i-1}]   \leq  \epsilon_n ^{\rho}\, \mathbb{E}[ \xi_{i}  ^{2}  | \mathcal{F}_{i-1}], \ \ \ 1\leq i \leq n.
\]
\item[(C2)]  There exists   $ \iota_n \in [0, \frac12]$ such that
$ \| \left\langle M\right\rangle _n-1\|_\infty \leq  \iota_n^2 .$
\end{description}
In many situations we have $ \epsilon_n, \iota_n\rightarrow 0$ as $n \rightarrow \infty$. In the case of sums of
i.i.d.\ random variables with finite $(2+\rho)$-th moments,  conditions (A1)  and (A2) are satisfied with $ \iota_n=0$ and $\epsilon_n=O(1/\sqrt{n} )$ as $n\rightarrow \infty.$

Define the self-normalized martingale
\begin{equation}
W_n = \frac{ M_n}{\sqrt{[ M ]_n}}, \quad n\geq 1.
\end{equation}
Define $\widehat{\epsilon}_m(x, \rho)$ in the same way as in (\ref{defepsilonm}) but with
$\varepsilon_m$ replaced by $\epsilon_m$.
   The proof of Theorem \ref{th1} is based on the  following technical lemma  which gives a Cram\'er type moderate deviation expansion for self-normalized martingales.

\begin{lemma}\label{lemma4.1}
Assume  conditions (C1) and (C2).
\begin{description}
  \item[\textbf{[i]}] If $\rho \in (0, 1)$, then   there exists an absolute constant $\alpha_{\rho,0} >0 $ such that for all $0\leq x \leq \alpha_{\rho,0 } \epsilon_n^{-1} $,
\begin{eqnarray*}
\Bigg|\ln \frac{\mathbb{P}( W_n  \geq x)}{1-\Phi \left( x\right)}  \Bigg|  \leq c_{\rho} \bigg( x^{2+\rho}  \epsilon_n^\rho+ x^2 \iota_n^2 +(1+x)\big(  \iota_n +\epsilon_n^\rho+ \widehat{\epsilon}_n(x, \rho) \big) \bigg).
\end{eqnarray*}

 \item[\textbf{[ii]}] If $\rho =1$, then   there exists an absolute constant $\alpha_0 >0 $ such that for all $0\leq x \leq \alpha_0  \epsilon_n^{-1} $,
\begin{eqnarray*}
\Bigg|\ln \frac{\mathbb{P}( W_n  \geq x)}{1-\Phi \left( x\right)}  \Bigg| \leq c  \bigg( x^{3}  \epsilon_n + x^2 \iota_n^2+(1+x)\big(   \iota_n+ \epsilon_n |\ln \epsilon_n| + \widehat{\epsilon}_n(x, 1)  \big) \bigg) .
\end{eqnarray*}

\end{description}
 Moreover,   the two above inequalities remain valid with $\frac{\mathbb{P}(W_n \leq -x)}{ \Phi \left( -x\right)}$
instead of  $\frac{\mathbb{P}( W_n \geq x)}{1-\Phi \left( x\right)}$.
\end{lemma}
\noindent\emph{Proof.} The points [i] and [ii] follows by  Corollary 2.3 of Fan \textit{et al.}\,\cite{FGLS17}.

\begin{remark}
  Notice that in Fan \emph{et al.}\ \cite{FGLS17}, the range for Lemma \ref{lemma4.1} is $0\leq x =o(  \epsilon_n^{-1})$. However,
  the proof of  Fan \emph{et al.}\ \cite{FGLS17} can be applied with no changes to extend
  the range to $0\leq x \leq \alpha_{\rho,0 }  \epsilon_n^{-1}$,
where $\alpha_{\rho,0 }$ is a sufficiently  small positive constant.
\end{remark}

Denote by $x^+=\max\{x, 0\}$  the positive part  of $x.$
\begin{lemma} \label{lemma4.2} Assume that  $\xi _i \geq -a$ a.s.\ for all $i\in [1, n]$. Write
\[
\emph{H}_n(\beta)  =\sum_{i=1}^n \Big(  \mathbb{E}\big[ (\xi_i^+)^\beta  |\mathcal{F}_{i-1} \big] +  a^\beta \Big), \ \ \ \ \beta \in (1,2]  .
\]
Then  for all $x, v>0$,
\begin{eqnarray}
 \mathbb{P}\left( S_n  \leq -x,\ \  \emph{H}_n(\beta) \leq v^\beta  \right)
 \leq   \exp\left\{- \frac{1}{2} C(\beta) \left(\frac{x}{v} \right)^\frac{\beta}{\beta -1} \right\},
\end{eqnarray}
where $C(\beta)= \beta^{-1 /(\beta-1) }  - \beta^{ -\beta /(\beta-1) } >0$ and $\beta \in (1, 2].$
\end{lemma}
\noindent\emph{Proof.}
Let $\beta \in (1, 2].$ Using the inequality
$$e^{-x} \leq 1- x +x^\beta \ \ \ \textrm{for} \ \ x \geq 0,$$
we have, for  all $i\in [1, n]$ and all $t>0, $
\begin{eqnarray*}
\mathbb{E}\big[  e^{-t  ( \xi_i + a ) }    \big| \mathcal{F}_{i-1}\big] & \leq& 1 - \mathbb{E}\big[  t  ( \xi_i + a )      \big| \mathcal{F}_{i-1}\big] + \mathbb{E}\big[  t^\beta  ( \xi_i + a )^{\beta}      \big| \mathcal{F}_{i-1}\big]  \\
&\leq& 1 - t a    +  2^{\beta-1} t^\beta \big(  \mathbb{E}\big[ (\xi_i^+)^\beta  |\mathcal{F}_{i-1} \big] +  a^\beta \big)\\
&\leq& \exp\{ - t a    +  2^{\beta-1} t^\beta \big(  \mathbb{E}\big[ (\xi_i^+)^\beta  |\mathcal{F}_{i-1} \big] +  a^\beta \big)  \}.
\end{eqnarray*}
Therefore, for all $x, t, v>0$,
\begin{eqnarray*}
&& \mathbb{P}\left( S_n  \leq -x,\ \  \textrm{H}_n(\beta) \leq v^\beta  \right) \\
&& \leq   \mathbb{E}\Big[ \exp\Big\{   - tx  -t \sum_{i=1}^n( \xi_i + a )  + t n a+ 2^{\beta-1}t^\beta \textrm{H}_n(\beta) -2^{\beta-1} t^\beta \textrm{H}_n(\beta) \Big \}  \mathbf{1}_{\{\textrm{H}_n(\beta) \leq v^\beta \}}  \Big] \\
&& \leq  e^{-tx + t n a+2^{\beta-1} t^\beta v^\beta } \mathbb{E}\Big[ \exp\Big\{   -t \sum_{i=1}^n( \xi_i + a )  - 2^{\beta-1}t^\beta \textrm{H}_n(\beta) \Big \} \Big] \\
&& \leq  e^{-tx + t n a+2^{\beta-1} t^\beta v^\beta } \mathbb{E}\Big[ \exp\Big\{   -t \sum_{i=1}^{n-1}( \xi_i + a )  -2^{\beta-1} t^\beta \textrm{H}_{n}(\beta) \Big \}  \mathbb{E}\big[  e^{-t  ( \xi_n + a ) }    \big| \mathcal{F}_{n-1}\big]\Big] \\
&& \leq  e^{-tx + t (n-1) a+ 2^{\beta-1}t^\beta v^\beta } \mathbb{E}\Big[ \exp\Big\{   -t \sum_{i=1}^{n-1}( \xi_i + a )  - 2^{\beta-1}t^\beta \textrm{H}_{n-1}(\beta) \Big \}  \Big] \\
&& \leq  e^{-tx +   2^{\beta-1}t^\beta v^\beta } .
\end{eqnarray*}
Taking $t=\frac{1 }{2} \big( \frac{x}{\beta v^\beta } \big)^{1/(\beta-1)} $
yields  the desired inequality.  \hfill\qed


The following  exponential inequality of  Peligrad et al.\ \cite{PUW07} (cf.\ Proposition 2 therein) plays an important role in the proof  of  Theorem \ref{th1}.
\begin{lemma}\label{lemma1}
  Let $(X_i)_{i \in \mathbb{Z}}$ be a sequence of random variables adapted to the filtration $(\mathcal{F}_i)_{i \in \mathbb{Z}}$.
  Then, for all $x \geq 0,$
  \begin{equation}
\mathbb{P}\bigg( \max_{1\leq i \leq n} |S_i| \geq x \bigg)  \leq 4 \sqrt{e} \exp  \Bigg\{ - \frac{x^2 }{ 2n (  \|X_1\|_\infty + 80 \sum_{j=1}^n j^{-3/2}\|\mathbb{E}[S_j|\mathcal{F}_0]\|_\infty )^2 } \Bigg\}.
\end{equation}
\end{lemma}

The last lemma shows that the tail probability of  $\max_{1\leq i \leq n} |S_i|$ has a sub-Gaussian decay rate.
 In the proof of Theorem  \ref{th1}, we apply it to estimate the tail probabilities for  the drift of a stationary sequence.

\subsection{Proof of Theorem  \ref{th1}}
Define
$$ D_{j}^\circ =S_{j}^\circ - \mathbb{E}[S_{j}^\circ | \mathcal{F}_{(j-1)m} ],\ \  1\leq j  \leq k.     $$
Then $(D_{j}^\circ, \mathcal{F}_{(j-1)m} )_{1\leq j \leq k}$ is a stationary  sequence of    martingale differences.
Clearly,
$$\mathbb{E}[ (D_j^\circ)^2 | \mathcal{F}_{(j-1)m}  ]=\mathbb{E}[ (S_j^\circ)^2 | \mathcal{F}_{(j-1)m}  ] -(\mathbb{E}[ S_j^\circ  | \mathcal{F}_{(j-1)m}  ])^2. $$
By stationarity  and the fact that $ k=\lfloor n/m   \rfloor,$   it follows that
$$ \frac1n \bigg\| \sum_{j=1}^k (\mathbb{E}[ S_j^\circ  | \mathcal{F}_{(j-1)m}  ])^2 \bigg\|_\infty \leq   \frac1m   \bigg\| \mathbb{E}[ S_{m}  | \mathcal{F}_{0}  ]\bigg\|_\infty^2  , $$
and that
\begin{eqnarray}
\bigg\| \frac1n \sum_{j=1}^k  \mathbb{E}[ (S_j^\circ)^2 | \mathcal{F}_{(j-1)m}  ] -   \sigma^2 \bigg\|_\infty
&\leq&   \frac1n \sum_{j=1}^k  \bigg\|  \mathbb{E}[ (S_j^\circ)^2 | \mathcal{F}_{(j-1)m}  ] -  m \sigma^2\bigg\|_\infty  + \frac{n-mk}n \sigma^2 \nonumber  \\
 &\leq&  \Big\|\frac1m \mathbb{E}[ S_m ^2 | \mathcal{F}_{0}  ] - \sigma^2\Big\|_\infty  +  \frac{m}n \sigma^2 .    \label{varsn}
\end{eqnarray}
Consequently, we have
\begin{eqnarray}
&& \bigg\| \frac1{n \sigma^2} \sum_{j=1}^{k }\mathbb{E}[ (D_j^\circ)^2 | \mathcal{F}_{(j-1)m}  ] - 1 \bigg\|_\infty \nonumber \\
&&\leq  \bigg\| \frac1{n \sigma^2} \sum_{j=1}^k  \mathbb{E}[ (S_j^\circ)^2 | \mathcal{F}_{(j-1)m}  ] -   1 \bigg\|_\infty     \ +\  \frac1{n \sigma^2} \bigg\| \sum_{j=1}^k (\mathbb{E}[ S_j^\circ  | \mathcal{F}_{(j-1)m}  ])^2 \bigg\|_\infty \nonumber \\
 && \leq \Big\|\frac1{m\sigma^2} \mathbb{E}[ S_m^2  | \mathcal{F}_{0}  ] -1\Big\|_\infty  +     \frac{m}n \ +\   \frac1{m\sigma^2}   \Big\| \mathbb{E}[ S_{m}  | \mathcal{F}_{0}  ]\Big\|_\infty^2 \nonumber \\
 && = \delta_m^2 +\frac{m}n  .  \label{thnsdf}
\end{eqnarray}
Since $\delta_m \rightarrow 0$ as $n \rightarrow \infty,$ it follows that
\begin{eqnarray}\label{fsfgk12}
 \Big\| \frac1m\mathbb{E}[ ( D_j^\circ)^{2} | \mathcal{F}_{(j-1)m}  ]  \Big\|_\infty  = \Big\|  \frac1m \mathbb{E}[ S_m^2  | \mathcal{F}_{0}  ] \textcolor{red}{-} \frac1m (\mathbb{E}[ S_m    | \mathcal{F}_{0}  ])^2  \Big\|_\infty  \sim \sigma^2,\ \ \ n\rightarrow \infty.
\end{eqnarray}
Using the inequality
\begin{eqnarray}
 |x-y|^{2+\rho}    \leq 2^{1+\rho} (|x|^{2+\rho}+|y|^{2+\rho}),  \label{jksineq}
 \end{eqnarray}
 by (\ref{fsfgk12}) and stationarity,  we deduce that
\begin{eqnarray}
  \mathbb{E}[ | D_j^\circ /(n^{1/2} \sigma) |^{2+\rho} | \mathcal{F}_{(j-1)m}  ]  &\leq & (n\sigma^2)^{-1- \rho/2}  2^{2+\rho}  \mathbb{E}[ |S_j^\circ|^{2+\rho} | \mathcal{F}_{(j-1)m}]  \nonumber \\
  &\leq&  \frac{2^{2+\rho}}{n^{\rho/2} \sigma^\rho   }   \bigg\|    \frac{  \mathbb{E}[ |S_j^\circ|^{2+\rho} | \mathcal{F}_{(j-1)m}] }{    \mathbb{E}[ ( D_j^\circ)^{2} | \mathcal{F}_{(j-1)m}  ]   }  \bigg\|_\infty  \mathbb{E}[ ( D_j^\circ/(n^{1/2} \sigma))^{2} | \mathcal{F}_{(j-1)m}  ]\nonumber \\
  &\leq& C_{\rho, 0} \frac{1}{n^{\rho/2}m \sigma^{2+\rho}   }   \Big\|   \mathbb{E}[ |S_m|^{2+\rho} | \mathcal{F}_{0}]   \Big\|_\infty  \mathbb{E}[ ( D_j^\circ/(n^{1/2} \sigma))^{2} | \mathcal{F}_{(j-1)m}  ]\nonumber \\
  &=  & C_{\rho, 0} \, \varepsilon_m^\rho \mathbb{E}[ ( D_j^\circ/(n^{1/2} \sigma))^{2} | \mathcal{F}_{(j-1)m}  ]. \label{sdfdf0}
\end{eqnarray}

 We first prove Theorem  \ref{th1} for $\rho \in (0, 1)$.  Set $\xi_j=  D_j^\circ /(n^{1/2} \sigma),  $  and denote $M_k=\sum_{j=1}^k \xi_j.$
 Then, by (\ref{thnsdf}) and (\ref{sdfdf0}), conditions (C1) and (C2) are satisfied with $n=k,$ $\epsilon_n^\rho =C_{\rho, 0} \varepsilon_m^\rho$
and   $\iota_n^2=\delta_m^2 +\frac{m}n$.
 By Lemma \ref{lemma4.1},
there exists a constant $\alpha_{\rho,0} >0 $ such that  for all $0\leq x \leq\alpha_{\rho,0}  \varepsilon_m^{-1} $,
\begin{eqnarray} \label{cram01}
&& \Bigg| \ln\frac{\mathbb{P}(M_k/\sqrt{[M]_k} \geq x)}{1-\Phi \left( x\right)}  \Bigg| \nonumber\\
  && \ \ \ \ \ \ \ \ \ \ \leq c_{\rho} \Bigg( x^{2+\rho}  \varepsilon_m^\rho+ x^2 (\delta_m^2  +\frac{m}n)+(1+x)\Big(  \delta_m  +\sqrt{\frac{m}n}+\varepsilon_m^{\rho} + \widehat{\varepsilon}_m(x, \rho) \Big) \Bigg)   .
\end{eqnarray}
Notice that,  by  Cauchy-Schwarz's inequality,
\begin{eqnarray*}
\bigg\| \frac{[M]_k}{    (V_k^\circ)^2 /(n \sigma^2)  }   -1  \bigg\|_\infty
  &=&  \bigg\| \frac{2}{ (V_k^\circ)^2}\sum_{j=1}^kS_{j}^\circ \mathbb{E}[ S_{j}^\circ | \mathcal{F}_{ (j-1)m } ]  +  \frac{1}{ (V_k^\circ)^2} \sum_{j=1}^k \Big(\mathbb{E}[S_{j}^\circ | \mathcal{F}_{(j-1)m} ] \Big)^2  \bigg\|_\infty \\
 &\leq&\bigg\|\frac{2}{(V_k^\circ)^2}\sum_{j=1}^kS_{j}^\circ \mathbb{E}[ S_{j}^\circ | \mathcal{F}_{ (j-1)m } ] \bigg\|_\infty +  \sum_{j=1}^k \bigg\|\frac{1}{(V_k^\circ)^2}(\mathbb{E}[S_{j}^\circ | \mathcal{F}_{(j-1)m} ])^2  \bigg\|_\infty \\
&\leq&     \bigg\|\frac{2}{ (V_k^\circ)^2 } \sum_{j=1}^k \big(\mathbb{E}[S_{j}^\circ | \mathcal{F}_{(j-1)m} ] \big)^2\bigg\|_\infty^{1/2} +  \sum_{j=1}^k \bigg\|\frac{1}{ V_k^\circ }\big|\mathbb{E}[S_{j}^\circ | \mathcal{F}_{(j-1)m} ]\big|  \bigg\|_\infty^2.
\end{eqnarray*}
 By stationarity and the fact that $\delta_m\rightarrow 0$,  when $  (V_k^\circ)^2   \geq \frac12 n \sigma^2,$ we have
\begin{eqnarray*}
 \bigg\| \frac{[M]_k}{    (V_k^\circ)^2 /(n \sigma^2)  }   -1  \bigg\|_\infty  &\leq& \frac{2 \sqrt{2}}{  \sqrt{n}  \sigma  } \bigg\| \sum_{j=1}^k \big(\mathbb{E}[S_{j}^\circ | \mathcal{F}_{(j-1)m} ] \big)^2\bigg\|_\infty^{1/2} +  \frac{2}{n \sigma^2   }\sum_{j=1}^k \Big\| \mathbb{E}[S_{j}^\circ | \mathcal{F}_{(j-1)m} ]   \Big\|_\infty^2  \\
   &\leq& \frac{2 \sqrt{2 k} }{  \sqrt{n}    \sigma  }  \Big\| \mathbb{E}[ S_m    | \mathcal{F}_{0}  ] \Big\|_\infty +  \frac{2}{ m \sigma^2}   \Big\| \mathbb{E}[ S_m    | \mathcal{F}_{0}  ] \Big\|_\infty^2 \\
 &\leq&    \frac 6{\sqrt{ m }\sigma }  \Big\| \mathbb{E}[ S_m    | \mathcal{F}_{0}  ] \Big\|_\infty=: \kappa_m.
\end{eqnarray*}
Clearly, $\delta_m \rightarrow 0 $   as $n\rightarrow\infty$  implies that $ \kappa_m \rightarrow 0$ as $n\rightarrow\infty.$
Thus, the last inequality implies that
\begin{eqnarray*}
 \frac{V_k^\circ }{\sqrt{n \sigma}}\geq \sqrt{\frac{[M]_k}{1+\kappa_m}} \geq\sqrt{[M]_k(1-\kappa_m)}   .
\end{eqnarray*}
 Recall that $ W^\circ_n= \frac{\sum_{j=1}^k S_{j}^\circ  }{ V_k^\circ   }=\frac{\sum_{j=1}^k S_{j}^\circ / (n^{1/2} \sigma) }{V_k^\circ/(n^{1/2} \sigma) }. $
It is easy to see that,  for all $x\geq 0,$
\begin{eqnarray}
 \mathbb{P}\bigg( W_n^\circ \geq x,   \, \, \frac{(V_k^\circ)^2}{n \sigma^2}    \geq \frac12  \bigg)    &\leq&   \mathbb{P}\bigg(\frac{\sum_{j=1}^k S_{j}^\circ / (n^{1/2} \sigma) }{\sqrt{[M]_k}}  \geq x \sqrt{1- \kappa_m} ,   \, \, \frac{(V_k^\circ)^2}{n \sigma^2}    \geq \frac12  \bigg) \nonumber  \\
  &\leq&   \mathbb{P}\bigg(\frac{M_k}{\sqrt{[M]_k}}  \geq x (1- \gamma_m |\ln \gamma_m| ) \sqrt{1- \kappa_m} ,   \, \, \frac{(V_k^\circ)^2}{n \sigma^2}    \geq \frac12  \bigg) \nonumber \\
  & &  + \mathbb{P}\bigg( \frac{1}{n^{1/2} \sigma}   \sum_{j=1}^k \mathbb{E}[S_{j}^\circ | \mathcal{F}_{(j-1)m} ]    \geq   x   \gamma_m |\ln \gamma_m|  \sqrt{1- \kappa_m}    \bigg) \nonumber  \\
  &\leq&   \mathbb{P}\bigg(\frac{M_k}{\sqrt{[M]_k}}  \geq x  (1- \gamma_m |\ln \gamma_m| )\sqrt{1- \kappa_m}   \bigg)  \nonumber \\
  & &  + \mathbb{P}\bigg(  \frac{1}{n^{1/2} \sigma}  \sum_{j=1}^k \mathbb{E}[S_{j}^\circ | \mathcal{F}_{(j-1)m} ]   \geq   x   \gamma_m |\ln \gamma_m|   \sqrt{1- \kappa_m}    \bigg)  \nonumber  \\
  & =:& I_1(x) + I_2(x).   \label{twomparts}
\end{eqnarray}

We proceed to estimate $I_1(x)$ and $I_2(x).$ First, we deal with $I_1(x).$
From (\ref{cram01}),  we have, for all $0\leq x \leq\alpha_{\rho,0}  \varepsilon_m^{-1} $,
\begin{eqnarray*}
 && \frac{I_1(x)}{1-\Phi \left( x (1- \gamma_m |\ln \gamma_m| )\sqrt{1-\kappa_m} \right)}  \nonumber  \\
   &&\ \ \ \  \ \ \ \  \ \  \ \ \  \leq \exp\bigg\{  c'_{\rho} \bigg( x^{2+\rho}  \varepsilon_m^\rho+ x^2 (\delta_m^2  +\frac{m}n)+(1+x)\Big(  \delta_m  +\sqrt{\frac{m}n} +\varepsilon_m^{\rho}+ \widehat{\varepsilon}_m(x, \rho) \Big) \bigg) \bigg\}.
\end{eqnarray*}
Using  the following inequalities
\begin{eqnarray}\label{fgsgj1}
\frac{1}{\sqrt{2 \pi}(1+x)} e^{-x^2/2} \leq 1-\Phi ( x ) \leq \frac{1}{\sqrt{ \pi}(1+x)} e^{-x^2/2}, \ \   x\geq 0,
\end{eqnarray}
we deduce  that,  for all $x\geq 0$ and $0\leq \varepsilon   \leq 1$,
\begin{eqnarray}
  \frac{1-\Phi \left( x \sqrt{1- \varepsilon} \right)}{1-\Phi \left( x\right) }& =& 1+ \frac{ \int_{x  \sqrt{1-  \varepsilon}  }^x \frac{1}{\sqrt{2\pi}}e^{-t^2/2}dt }{1-\Phi \left( x\right) } \nonumber \\
  &\leq&    1+ \frac{\frac{1}{\sqrt{2\pi}} e^{-x^2(1- \varepsilon)/2}  x \varepsilon  }{ \frac{1}{\sqrt{2 \pi} (1+x)} e^{-x^2/2}  }  \nonumber \\
   &\leq &   1+  C (1+  x^2)    \varepsilon  e^{     x^2 \varepsilon /2 }   \nonumber \\
    &\leq&     \exp\Big\{ C  (1+  x^2 )  \varepsilon  \Big\}. \label{sfdsh}
\end{eqnarray}
Notice that  $(1- \gamma_m |\ln \gamma_m| )\sqrt{1- \kappa_m} \geq \sqrt{1-2(\gamma_m |\ln \gamma_m|+\kappa_m)}. $
Using   inequality (\ref{sfdsh}) and  the fact that  $  \kappa_m\leq 6  \gamma_m \leq 6  \gamma_m|\ln \gamma_m|$, we obtain, for all $0\leq x \leq\alpha_{\rho,0}  \varepsilon_m^{-1} $,
\begin{eqnarray}
 && \frac{I_1(x)}{1-\Phi \left( x \right)} =  \frac{I_1(x)}{1-\Phi \left( x (1- \gamma_m |\ln \gamma_m| )\sqrt{1- \kappa_m} \right)}    \frac{1-\Phi \left( x (1- \gamma_m |\ln \gamma_m| )\sqrt{1-\kappa_m} \right)}{1-\Phi \left( x \right)}  \nonumber   \\
  && \ \ \   \leq  \exp\Bigg\{  C'_{\rho} \bigg( x^{2+\rho}  \varepsilon_m^\rho+ x^2 \Big(\delta_m^2+\frac{m}n +\gamma_m |\ln \gamma_m| + \kappa_m\Big)  \nonumber \\
  &&  \ \ \   \ \ \   \ \ \   \ \ \  \ \ \   \ \ \   \ \ \      + (1+x)\Big(  \delta_m +\sqrt{\frac{m}n}+\varepsilon_m^{\rho} + \widehat{\varepsilon}_m(x, \rho)+\gamma_m |\ln \gamma_m|+\kappa_m \Big) \bigg) \Bigg\} \nonumber  \\
   && \ \ \  \leq  \exp\Bigg\{  C''_{\rho} \bigg( x^{2+\rho}  \varepsilon_m^\rho+ x^2 \Big(\delta_m^2   +\gamma_m |\ln \gamma_m|+\frac{m}n \Big)    \nonumber \\
  &&  \ \ \   \ \ \   \ \ \   \ \ \  \ \ \   \ \ \   \ \ \   \ \ \   \ \ \  \ \ \    + (1+x)\Big(  \delta_m+\sqrt{\frac{m}n}+ \varepsilon_m^{\rho}+\gamma_m |\ln \gamma_m| + \widehat{\varepsilon}_m(x, \rho) \Big) \bigg) \Bigg\},  \label{gfdsg}
\end{eqnarray}
which gives the suitable bound for $I_1(x).$

Now  we deal with $I_2(x).$
By Lemma \ref{lemma1},  the definition of $\gamma_m$ (cf.\ (\ref{grmman})) and the fact that $\gamma_m \rightarrow 0$,    we derive that, for all $x \geq 0$,
\begin{eqnarray}
I_2(x) & \leq&
 4 \sqrt{e} \exp  \bigg\{ - \frac{n \sigma^2 x^2\gamma_m^2 (\ln \gamma_m)^2 (1-\kappa_m) }{ 2 k ( \big\| \mathbb{E}[ S_{m}  | \mathcal{F}_{0}  ]\big\|_\infty + 80 \sum_{j=1}^k j^{-3/2}\|\mathbb{E}[S_{jm}|\mathcal{F}_0]\|_\infty )^2 } \bigg\}  \nonumber\\
  &\leq&
 4 \sqrt{e} \exp  \bigg\{ - C_0 x^2 (\ln \gamma_m)^2  \bigg\}.  \label{ineq6}
\end{eqnarray}
From the last inequality,  using  (\ref{fgsgj1}),
  we deduce that, for all $x \geq 1,$
\begin{eqnarray}
   \frac{I_2(x)}{1-\Phi \left( x \right)} & \leq& C_1 (1+x)  \exp  \bigg\{ - C_0 x^2 (\ln \gamma_m)^2 + \frac12 x^2  \bigg\} \label{fsdfds312}\\
    \ \ \ \ \  &\leq& C_2 \, (1+x)   \gamma_m |\ln \gamma_m| , \nonumber
\end{eqnarray}
which gives the suitable bound for $I_2(x).$
Thus,  from (\ref{twomparts}), for all   $x \geq 1,$
\begin{eqnarray}
 && \frac{ \mathbb{P}\Big( W_n^\circ \geq x,\, \frac{(V_k^\circ)^2}{n \sigma^2}  \geq \frac12 \Big)}{1-\Phi \left( x \right)}  \leq   \frac{  I_1(x) + I_2(x)   }{1-\Phi \left( x \right)} \nonumber \\
  &&   \leq  \exp\Bigg\{  c''_{\rho} \bigg( x^{2+\rho}  \varepsilon_m^\rho+ x^2 \Big(\delta_m^2   +\gamma_m |\ln \gamma_m|+\frac{m}n \Big) \nonumber \\
    &&\ \ \ \ \ \ \ \ \  \ \ \ \ \ \ \ \ \ \ \  \ \ \ \ \ \   + (1+x)\Big(  \delta_m  +\sqrt{\frac{m}n}+\varepsilon_m^{\rho}+\gamma_m |\ln \gamma_m| + \widehat{\varepsilon}_m(x, \rho) \Big) \bigg) \Bigg\} .  \label{dgsxcg}
\end{eqnarray}

Clearly,   we have
\begin{eqnarray}
\mathbb{P}\bigg( (V_k^\circ)^2     < \frac12   n \sigma^2   \bigg) & =& \mathbb{P}\Bigg(   \sum_{j=1}^k \Big((S_{j}^\circ)^2 -  \mathbb{E}[ (S_j^\circ)^2 | \mathcal{F}_{(j-1)m}  ] \Big)  < \frac12   n \sigma^2   -\sum_{j=1}^k \mathbb{E}[ (S_j^\circ)^2 | \mathcal{F}_{(j-1)m}  ]  \Bigg) \nonumber \\
&   \leq& \mathbb{P}\Bigg(   \sum_{j=1}^k \Big((S_{j}^\circ)^2 -  \mathbb{E}[ (S_j^\circ)^2 | \mathcal{F}_{(j-1)m}  ] \Big)  < -\frac14   n \sigma^2    \Bigg) ,\label{ineq30d}
\end{eqnarray}
where the last line follows by (\ref{varsn})  and the fact that $\delta_m\rightarrow 0$ and $m/n \rightarrow0$.
Denote $$\eta_j=\Big(\frac{S_{j}^\circ}{\sigma \sqrt{n} } \Big)^2 -  \mathbb{E}\Big[ \Big(\frac{S_{j}^\circ}{\sigma \sqrt{n} } \Big)^2 \Big| \mathcal{F}_{(j-1)m}  \Big]. $$
Then,   by (\ref{jksineq}) and stationarity,  it is easy to see that
\begin{eqnarray*}
\sum_{j=1}^k\Big\|\mathbb{E}[|\eta_j|^{(2+\rho)/2} |\mathcal{F}_{(j-1)m}] \Big\|_\infty  \leq 2^{1+\rho} \sum_{j=1}^k\Big\|  \mathbb{E}\Big[ \Big|\frac{S_{j}^\circ}{\sigma \sqrt{n} } \Big|^{2+\rho} \Big| \mathcal{F}_{(j-1)m}  \Big]\Big\|_\infty  \leq  2^{1+\rho}  \varepsilon_m^\rho
\end{eqnarray*}
and that, for some positive constant $c,$
$$   \eta_i \geq - \frac{1}{n\sigma^2} \big\|\mathbb{E}\big[  S_m^2 \big| \mathcal{F}_{0}  \big] \big\|_\infty\geq -\frac mn c \ \ \ \ \textrm{a.s.},$$
where the last inequality follows by the fact that $\delta_m \rightarrow 0$ as $n\rightarrow\infty.$
From (\ref{ineq30d}), using Lemma  \ref{lemma4.2} with $a= \frac mn c$ and $\beta=(2+\rho)/2$, we have
\begin{eqnarray}
\mathbb{P}\bigg(   \frac{(V_k^\circ)^2}{n \sigma^2}  < \frac12  \bigg)
\leq   \exp\bigg\{- C(\rho) \Big( \frac {1}{ \varepsilon_m^{2}}  + \frac n m \Big)  \bigg\}, \label{dfsf023}
\end{eqnarray}
where $C(\rho) >0$ depends only on $\rho. $  Notice that, by (\ref{fgsgj1}), it holds, for small enough $\alpha_{\rho,0}>0$ and all $1\leq x \leq  \alpha_{\rho,0} \min\{\varepsilon_n^{-1},\, \sqrt{n/m} \} $,
\begin{eqnarray}
\frac{ 1}{1-\Phi \left( x \right)} \exp\Bigg\{- C(\rho) \Big( \frac {1}{ \varepsilon_m^{2}}  + \frac n m \Big)  \Bigg\}
 \leq   \sqrt{2 \pi}(1+x)\Big(  \sqrt{\frac{m}n}+\varepsilon_m^{\rho}   \Big). \label{dfssf0}
\end{eqnarray}
Then,  by (\ref{dgsxcg}), (\ref{dfsf023}) and (\ref{dfssf0}),  we obtain,  for all $1\leq x \leq  \alpha_{\rho,0} \min\{\varepsilon_n^{-1},\, \sqrt{n/m} \}  , $
\begin{eqnarray*}
 && \frac{ \mathbb{P}\Big( W_n^\circ \geq x \Big)}{1-\Phi \left( x \right)}  \leq  \frac{ \mathbb{P}\Big( W_n^\circ \geq x,\, \frac{(V_k^\circ)^2}{n \sigma^2}  \geq \frac12 \Big)}{1-\Phi \left( x \right)} + \frac{ \mathbb{P}\Big( \frac{(V_k^\circ)^2}{n \sigma^2}  < \frac12 \Big)}{1-\Phi \left( x \right)} \\
  &&   \leq  \exp\Bigg\{  c''_{\rho} \bigg( x^{2+\rho}  \varepsilon_m^\rho+ x^2 \Big(\delta_m^2   +\gamma_m |\ln \gamma_m|+\frac{m}n \Big)  \\
    && \quad \quad  \quad   \quad \quad  \quad \quad  \quad \quad + (1+x)\Big(  \delta_m +\sqrt{\frac{m}n} +\varepsilon_m^{\rho}+\gamma_m |\ln \gamma_m| + \widehat{\varepsilon}_m(x, \rho) \Big) \bigg) \Bigg\}
    \\
    && \ \ \ \  + \ \frac{ 1}{1-\Phi \left( x \right)} \exp\Bigg\{- C(\rho) \Big( \frac {1}{ \varepsilon_m^{2}}  + \frac n m \Big)  \Bigg\} \\
  &&   \leq  \exp\Bigg\{  c'''_{\rho} \bigg( x^{2+\rho}  \varepsilon_m^\rho+ x^2 \Big(\delta_m^2   +\gamma_m |\ln \gamma_m|+\frac{m}n \Big)  \\
    && \quad \quad  \quad  \quad \quad  \quad \quad  \quad \quad + (1+x)\Big(  \delta_m +\sqrt{\frac{m}n}+\varepsilon_m^{\rho}+\gamma_m |\ln \gamma_m|  + \widehat{\varepsilon}_m(x, \rho) \Big) \bigg) \Bigg\}.
\end{eqnarray*}
From the last inequality, we get,
for all $1\leq x \leq \alpha_{\rho   } \min\{\varepsilon_m^{-1},\, \sqrt{n/m} \}  , $
\begin{eqnarray}
  \ln \frac{\mathbb{P}( W_n^\circ  \geq x   )}{1-\Phi \left(  x\right)}     &\leq&  c'''_{\rho} \bigg( x^{2+\rho}  \varepsilon_m^\rho+ x^2 \Big(\delta_m^2   +\gamma_m |\ln \gamma_m|+\frac{m}n \Big)   \nonumber \\
 &&\ \ \ \ \ \  \  +\, (1+x)\Big(  \delta_m +\sqrt{\frac{m}n}+\varepsilon_m^{\rho}+\gamma_m |\ln \gamma_m|  + \widehat{\varepsilon}_m(x, \rho) \Big) \bigg)  , \label{fgghs3}
\end{eqnarray}
which gives the upper bound of $\ln \frac{\mathbb{P}( W_n^\circ  \geq x   )}{1-\Phi \left(  x\right)}$ for $\rho \in (0, 1)$.
The proof of the lower bound of $\ln \frac{\mathbb{P}( W_n^\circ  \geq x   )}{1-\Phi \left(  x\right)}$,  $1\leq x \leq  \alpha_{\rho,0} \min\{\varepsilon_n^{-1},\, \sqrt{n/m} \} $, is  similar to the proof of (\ref{fgghs3}), but,  instead of using (\ref{twomparts}), we 
use the following inequalities:  for all $x\geq 0,$
\begin{eqnarray}
 \mathbb{P}\bigg( W_n^\circ \geq x  \bigg) &\geq& \mathbb{P}\bigg( W_n^\circ \geq x,   \, \, \frac{(V_k^\circ)^2}{n \sigma^2}    \geq \frac12  \bigg) \nonumber \\
   &\geq&   \mathbb{P}\bigg(\frac{\sum_{j=1}^k S_{j}^\circ / (n^{1/2} \sigma) }{\sqrt{[M]_k}}  \geq x \sqrt{1+\kappa_m} ,   \, \, \frac{(V_k^\circ)^2}{n \sigma^2}    \geq \frac12  \bigg) \nonumber\\
  &\geq&   \mathbb{P}\bigg(\frac{M_k}{\sqrt{[M]_k}}  \geq x (1+ \gamma_m |\ln \gamma_m| ) \sqrt{1+\kappa_m} ,   \, \, \frac{(V_k^\circ)^2}{n \sigma^2}    \geq \frac12  \bigg) \nonumber\\
  & &  - \ \mathbb{P}\bigg( \frac{1}{n^{1/2} \sigma} \sum_{j=1}^k \mathbb{E}[S_{j}^\circ | \mathcal{F}_{(j-1)m} ]    \geq   x   \gamma_m |\ln \gamma_m|   \sqrt{1+\kappa_m},      \, \, \frac{(V_k^\circ)^2}{n \sigma^2}    \geq \frac12     \bigg)\nonumber\\
   &\geq&   \mathbb{P}\bigg(\frac{M_k}{\sqrt{[M]_k}}  \geq x (1+ \gamma_m |\ln \gamma_m| ) \sqrt{1+\kappa_m}   \bigg)- \mathbb{P}\bigg( \frac{(V_k^\circ)^2}{n \sigma^2}    < \frac12  \bigg) \nonumber\\
  & &- \ \mathbb{P}\bigg( \frac{1}{n^{1/2} \sigma} \sum_{j=1}^k \mathbb{E}[S_{j}^\circ | \mathcal{F}_{(j-1)m} ]    \geq   x   \gamma_m |\ln \gamma_m|   \sqrt{1+\kappa_m}    \bigg) \nonumber \\
  &=:&P_1(x)-P_2-P_3(x). \label{fshlmo}
\end{eqnarray}
By an argument similar to that of (\ref{gfdsg}), we deduce that, for all $0\leq x \leq\alpha_{\rho,0}  \varepsilon_m^{-1} $,
\begin{eqnarray}
&& \frac{P_1(x) }{1-\Phi(x)}   \geq \exp\bigg\{-  c_{\rho} \bigg( x^{2+\rho}  \varepsilon_m^\rho+ x^2 (\delta_m^2+ \gamma_m |\ln \gamma_m|+\frac{m}n)   \nonumber \\
&&\ \ \ \ \ \ \ \ \ \ \ \ \ \ \ \ \ \ \  \ \ \ \ \ \ \   +(1+x)\Big(  \delta_m  +\sqrt{\frac{m}n}+\varepsilon_m^{\rho} + \gamma_m |\ln \gamma_m|+ \widehat{\varepsilon}_m(x, \rho) \Big) \bigg) \bigg\}  .
\end{eqnarray}
By (\ref{dfsf023}), we have, for small enough $\alpha_{\rho,0}>0$ and all $0\leq x \leq  \alpha_{\rho,0} \min\{\varepsilon_n^{-1},\, \sqrt{n/m} \} $,
\begin{eqnarray}
\frac{P_2}{1-\Phi \left( x \right)} & \leq&  \sqrt{2 \pi} (1+x)     \exp\bigg\{- C(\rho) \Big( \frac {1}{ \varepsilon_m^{2}}  + \frac n m \Big) +\frac12 x^2  \bigg\} \nonumber \\
&\leq& C_{\rho,3}\Big(  \sqrt{\frac{m}n}+\varepsilon_m^{\rho}   \Big)\exp  \bigg\{ -  \frac12 x^2  \bigg\}.
\end{eqnarray}
By an argument similar to that of (\ref{fsdfds312}), we get, for all $x \geq 1,$
\begin{eqnarray}
   \frac{P_3(x)}{1-\Phi \left( x \right)} & \leq& C_1 (1+x)  \exp  \bigg\{ - C_0 x^2 (\ln \gamma_m)^2 + \frac12 x^2  \bigg\} \nonumber \\
    \ \ \ \ \  &\leq& C_4 \,    \gamma_m |\ln \gamma_m| \exp  \bigg\{ -  \frac12 x^2  \bigg\}  . \label{fqsdfs2}
\end{eqnarray}
Combining the inequalities (\ref{fshlmo})-(\ref{fqsdfs2}) together,  we obtain, for all $1\leq x \leq  \alpha_{\rho,0} \min\{\varepsilon_n^{-1},\, \sqrt{n/m} \} $,
\begin{eqnarray*}
  \ln \frac{\mathbb{P}( W_n^\circ  \geq x   )}{1-\Phi \left(  x\right)}     &\geq&  - C_{\rho} \bigg( x^{2+\rho}  \varepsilon_m^\rho+ x^2 \Big(\delta_m^2   +\gamma_m |\ln \gamma_m|+\frac{m}n \Big)   \nonumber \\
 &&\ \ \ \ \ \  \ \ \ \ \ +\, (1+x)\Big(  \delta_m +\sqrt{\frac{m}n}+\varepsilon_m^{\rho}+\gamma_m |\ln \gamma_m|  + \widehat{\varepsilon}_m(x, \rho) \Big) \bigg) .
\end{eqnarray*}
This completes the proof of Theorem  \ref{th1} for all $1\leq x \leq  \alpha_{\rho,0} \min\{\varepsilon_n^{-1},\, \sqrt{n/m} \} $.

 For  the case  $0\leq x \leq 1,$   instead of  (\ref{twomparts}), we make use of the following estimations:
\begin{eqnarray*}
 \mathbb{P}\bigg( W_n^\circ \geq x,   \, \, \frac{(V_k^\circ)^2}{n \sigma^2}    \geq \frac12  \bigg)    &\leq&   \mathbb{P}\bigg(\frac{\sum_{j=1}^k S_{j}^\circ / (n^{1/2} \sigma) }{\sqrt{[M]_k}}  \geq x \sqrt{1- \kappa_m} ,   \, \, \frac{(V_k^\circ)^2}{n \sigma^2}    \geq \frac12  \bigg) \\
  &\leq&   \mathbb{P}\bigg(\frac{M_k}{\sqrt{[M]_k}}  \geq (  x - \gamma_m |\ln \gamma_m| ) \sqrt{1- \kappa_m} ,   \, \, \frac{(V_k^\circ)^2}{n \sigma^2}    \geq \frac12  \bigg)   +  I_2(1) \\
  &\leq&   \mathbb{P}\bigg(\frac{M_k}{\sqrt{[M]_k}}  \geq    (x- \gamma_m |\ln \gamma_m| )\sqrt{1- \kappa_m}   \bigg) +  I_2(1) \\
  & =:& \widetilde{I}_1(x) +  I_2(1).
\end{eqnarray*}
By an argument  similar to the case of  $1\leq x \leq  \alpha_{\rho,0} \min\{\varepsilon_n^{-1},\, \sqrt{n/m} \} $, we obtain the upper bound of
$\displaystyle \ln \frac{\mathbb{P}( W_n^\circ  \geq x   )}{1-\Phi \left(  x\right)}$ for all $0\leq x \leq 1$.
To prove the lower bound of $\displaystyle \ln \frac{\mathbb{P}( W_n^\circ  \geq x   )}{1-\Phi \left(  x\right)}, 0\leq x \leq 1,$ instead of  (\ref{twomparts}),  we should  use  the following estimations:
\begin{eqnarray*}
\mathbb{P}\bigg( W_n^\circ \geq x   \bigg)  &\geq&
 \mathbb{P}\bigg( W_n^\circ \geq x,   \, \, \frac{(V_k^\circ)^2}{n \sigma^2}    \geq \frac12  \bigg)  \\
   &\geq&   \mathbb{P}\bigg(\frac{\sum_{j=1}^k S_{j}^\circ / (n^{1/2} \sigma) }{\sqrt{[M]_k}}  \geq x \sqrt{1+ \kappa_m} ,   \, \, \frac{(V_k^\circ)^2}{n \sigma^2}    \geq \frac12  \bigg) \\
  &\geq&   \mathbb{P}\bigg(\frac{M_k}{\sqrt{[M]_k}}  \geq (  x + \gamma_m |\ln \gamma_m| ) \sqrt{1+\kappa_m} ,   \, \, \frac{(V_k^\circ)^2}{n \sigma^2}    \geq \frac12  \bigg)    -   I_2(1) \\
  & \geq& \mathbb{P}\bigg(\frac{M_k}{\sqrt{[M]_k}}  \geq (  x + \gamma_m |\ln \gamma_m| ) \sqrt{1+\kappa_m}    \bigg) - \mathbb{P}\bigg(  \frac{(V_k^\circ)^2}{n \sigma^2}   <\frac12  \bigg)  -   I_2(1).
\end{eqnarray*}
Again by an argument  similar to the case of  $1\leq x \leq  \alpha_{\rho,0} \min\{\varepsilon_n^{-1},\, \sqrt{n/m} \} $, we get the lower bound of
$\displaystyle \ln \frac{\mathbb{P}( W_n^\circ  \geq x   )}{1-\Phi \left(  x\right)}$ for all $0\leq x \leq 1$.
This completes the proof of Theorem  \ref{th1}  for   $\rho \in (0, 1)$.

 For $\rho=1,$ the proof of Theorem  \ref{th1} is similar to the case of $\rho \in (0, 1)$, where the term $\varepsilon_m |\ln \varepsilon_m|$ comes from point [ii] of Lemma \ref{lemma4.1} with $\epsilon_n  =C_{1, 0} \varepsilon_m.$
Notice that if $(X_i)_{i \in \mathbb{Z}}$  satisfies the condition   of   Theorem  \ref{th1}, then
  $(-X_i)_{i \in \mathbb{Z}}$ also satisfies the same condition.
   Thus the assertions in  Theorem  \ref{th1}
 remain valid when $\displaystyle \frac{\mathbb{P}( W_n^\circ  \geq x  )}{1-\Phi \left(  x\right)}$  is replaced by $\displaystyle \frac{\mathbb{P}(W_n^\circ \leq -x )}{ \Phi \left(-  x\right)}$, $x\geq 0.$

\subsection{Proof of Corollary \ref{co0} }\label{sec4.3}
First, we prove that
\begin{eqnarray}\label{dfgk18}
 \limsup_{n\rightarrow \infty} a_n^2 \ln \mathbb{P}\bigg(  a_n W^\circ_n\in B  \bigg) \leq  - \inf_{x \in \overline{B}}\frac{x^2}{2   }.
\end{eqnarray}
For any given Borel set $B\subset \mathbb{R},$ let $x_0=\inf_{x\in B} |x| \geq\inf_{x\in \overline{B}} |x|.$
By Theorem \ref{th1}, we deduce that
\begin{eqnarray*}
 && \mathbb{P}\bigg(  a_n W^\circ_n \in B \bigg) \\
 &&\leq  \mathbb{P}\bigg(\, W^\circ_n \geq  \frac{ x_0}{a_n }\bigg)+  \mathbb{P}\bigg(\, W^\circ_n \leq-  \frac{ x_0}{a_n }\bigg)\\
 &&\leq  2\bigg( 1-\Phi \Big( \frac{ x_0}{a_n  }\Big)\bigg) \exp\Bigg\{C \bigg( (\frac{ x_0}{a_n  })^{2+\rho} \varepsilon_m^\rho + \Big (\frac{ x_0}{a_n  } \Big)^2 \Big(\delta_m^2 + \gamma_m|\ln \gamma_m| +  \frac m n\Big) \\
  && \  \  \ \ \ \    \  \  \ \ \ \  \ \ \ \ \ \   \  \  \ \ \ \    \  \  \ \ \ \    \  \  \  \  \  \ \ \ \    \   +\, (1+ \frac{ x_0}{a_n })\Big(  \delta_m  +\gamma_m|\ln \gamma_m|+\varepsilon_m^{\rho/4}+ \sqrt{\frac m n}\Big)\bigg) \Bigg\}.
\end{eqnarray*}
Notice that  $a_n\rightarrow 0$ and $ a_n  \min\{\varepsilon_m^{-1},\, \sqrt{n/m} \}\rightarrow \infty$ as $n\rightarrow \infty.$
Using   (\ref{fgsgj1}) and (\ref{bhg03}),
we deduce that
\begin{eqnarray*}
\limsup_{n\rightarrow \infty} a_n^2 \ln \mathbb{P}\bigg( a_n  W^\circ_n  \in B  \bigg)
 \ \leq \  -\frac{x_0^2}{2 } \ \leq \  - \inf_{x \in \overline{B}}\frac{x^2}{2 } ,
\end{eqnarray*}
which gives (\ref{dfgk18}).

Next, we prove that
\begin{eqnarray}\label{dfgk02}
\liminf_{n\rightarrow \infty} a_n^2 \ln \mathbb{P}\bigg(  a_n  W^\circ_n \in B  \bigg) \geq   - \inf_{x \in B^o}\frac{x^2}{ 2   } .
\end{eqnarray}
Without loss of generality, we assume that $B^o \neq \emptyset,$ otherwise (\ref{dfgk02})  holds obviously, \textcolor{red}{since} in this case the infimum of a function over an empty set is
equal to $\infty$  by convention.
For any given $\varepsilon_1>0,$ there exists an $x_0 \in B^o $ such that
\begin{eqnarray}
 0< \frac{x_0^2}{2 } \leq   \inf_{x \in B^o}\frac{x^2}{2  } +\varepsilon_1.
\end{eqnarray}
We only consider the case when $x_0>0$, the case $x_0<0$ being proved in the same way. Since $B^o$ is an open set,
for $x_0 \in B^o$ and  small enough  $\varepsilon_2 \in (0, x_0) $ it holds $(x_0-\varepsilon_2, x_0+\varepsilon_2]  \subset B.$
Clearly, $x_0\geq\inf_{x\in \overline{B}} x.$
It is easy to see  that
\begin{eqnarray*}
\mathbb{P}\bigg(  a_n  W^\circ_n \in B  \bigg)   &\geq&   \mathbb{P}\bigg(   W^\circ_n  \in ( a_n^{-1}  ( x_0-\varepsilon_2), a_n^{-1} ( x_0+\varepsilon_2)] \bigg)\\
&\geq&   \mathbb{P}\Big(   W^\circ_n   > a_n^{-1}( x_0-\varepsilon_2)   \Big)-\mathbb{P}\Big(  W^\circ_n >   a_n^{-1}( x_0+\varepsilon_2) \Big).
\end{eqnarray*}
 By Theorem \ref{th1}, we have
 $$\lim_{n\rightarrow \infty} \frac{\mathbb{P}\Big( W_n^o > a_n^{-1}( x_0+\varepsilon_2) \Big)}{\mathbb{P}\Big(  W_n ^o  >  a_n^{-1} ( x_0-\varepsilon_2)   \Big) } =0 .$$
Again, by Theorem \ref{th1},  (\ref{fgsgj1}) and (\ref{bhg03}), it follows that
\begin{eqnarray*}
\liminf_{n\rightarrow \infty} a_n^2\ln \mathbb{P}\bigg(a_n  W^\circ_n \in B  \bigg) \geq   \liminf_{n\rightarrow \infty} a_n^2\ln \frac12\mathbb{P}\Big(   W^\circ_n   > a_n^{-1}( x_0-\varepsilon_2)   \Big)     = -  \frac{1}{2 }( x_0-\varepsilon_2)^2 . \label{ffhms}
\end{eqnarray*}
Letting $\varepsilon_2\rightarrow 0,$  we obtain
\begin{eqnarray*}
\liminf_{n\rightarrow \infty}a_n^2\ln \mathbb{P}\bigg(a_n  W^\circ_n   \in B \bigg) &\geq& - \frac{x_0^2}{2 }  \  \geq \   -\inf_{x \in B^o}\frac{x^2}{2 } -\varepsilon_1.
\end{eqnarray*}
Since $\varepsilon_1>0$ can be arbitrarily small, we get (\ref{dfgk02}).  The proof of Corollary \ref{co0} is complete.

\subsection{Proof of Corollary \ref{co01} }\label{sec4.4}
We only need to consider the case where $\max\{\gamma_m,  \varepsilon_m,  \delta_m, m/n \}  \leq  1/10.$ Otherwise,
Corollary \ref{co01} holds obviously by choosing $C_\rho$ large enough.
Denote $$\kappa_n=  \min\{\gamma_m^{-1/4} ,\, \varepsilon_m^{-\rho(2-\rho)/8},  \delta_m^{-1/4}, (m/n)^{-1/4}   \} .  $$
It is easy to see that
\begin{eqnarray}
\sup_{ x   }  \Big|\mathbb{P}( W^\circ_n   \leq x   )  -  \Phi \left( x\right) \Big| &\leq&  \sup_{  |x|\leq \kappa_n } \Big|\mathbb{P}(W^\circ_n   \leq x  )  -  \Phi \left( x\right) \Big|
   + \sup_{  |x|> \kappa_n} \Big|\mathbb{P}( W^\circ_n   \leq x  )  -  \Phi \left( x\right) \Big|  \nonumber\\
&= & \sup_{  |x|\leq \kappa_n } \Big|\mathbb{P}( W^\circ_n   \leq x  )  -  \Phi \left( x\right) \Big|  \nonumber\\
&&  + \sup_{   x < - \kappa_n}  \mathbb{P}( W^\circ_n   \leq x  )  +  \sup_{   x < - \kappa_n}  \Phi \left( x\right)    \nonumber\\
&&  + \sup_{   x >   \kappa_n}  \mathbb{P}( W^\circ_n   >x  )  +  \sup_{   x > \kappa_n} ( 1- \Phi \left( x\right) ).   \label{ineq010}
\end{eqnarray}
 Notice that $$\sup_{ |x|\leq \kappa_n }\{\varepsilon_m ^\rho|\ln \varepsilon_m|, \widehat{\varepsilon}_m(x, \rho)\}  =\varepsilon_m^{ \rho(2-\rho)/4 }.$$
By Theorem \ref{th1} and the inequality $|e^x-1|\leq |x|e^{|x|},$ we have
\begin{eqnarray}
 &&\sup_{ |x|\leq \kappa_n } \Big|\mathbb{P}( W^\circ_n   \leq x  )  -  \Phi \left( x\right) \Big| \nonumber \\
  &&\leq  \sup_{ |x|\leq \kappa_n } \Big(1-\Phi(|x|) \Big) \bigg| e^{   C_{\rho}  \big( x^{2+\rho}  \varepsilon_m^\rho+ x^2 \big(\delta_m^2   +\gamma_m |\ln \gamma_m|+\frac{m}n \big)    \  +\ (1+x)\big(  \delta_m +\gamma_m |\ln \gamma_m|  +\varepsilon_m^{ \rho(2-\rho)/4 }+\sqrt{\frac{m}n} \big)\big)}   -1 \bigg| \nonumber\\
 &&\leq  C_{\rho, 1 }  \Big( \delta_m +\gamma_m |\ln \gamma_m| +\varepsilon_m^{ \rho(2-\rho)/4 }+\sqrt{\frac{m}n}   \Big ). \label{ineq020}
 \end{eqnarray}
From the last inequality, we get
 \begin{eqnarray}
\sup_{   x < - \kappa_n}  \mathbb{P}( W^\circ_n   \leq x   )&= &  \mathbb{P}( W^\circ_n   \leq  - \kappa_n   ) \nonumber  \\
&\leq& C_{\rho, 1 } \Big( \delta_m +\gamma_m |\ln \gamma_m| +\varepsilon_m^{ \rho(2-\rho)/4 }+\sqrt{\frac{m}n}   \Big ) +  \Phi \left( - \kappa_n\right) \nonumber \\
&\leq& C_{\rho, 2 } \Big(  \delta_m +\gamma_m |\ln \gamma_m| +\varepsilon_m^{ \rho(2-\rho)/4 }+\sqrt{\frac{m}n}   \Big ).   \label{ineq030}
\end{eqnarray}
Similarly, we have
\begin{eqnarray}   \label{ineq040}
\sup_{   x >   \kappa_n}  \mathbb{P}(W^\circ_n     > x   )    &  \leq&  C_{\rho, 3 } \Big(  \delta_m +\gamma_m |\ln \gamma_m| +\varepsilon_m^{ \rho(2-\rho)/4 }+\sqrt{\frac{m}n}   \Big ).
\end{eqnarray}
Clearly, it holds that
\begin{eqnarray} \label{ineq050}
  \sup_{   x > \kappa_n} ( 1- \Phi \left( x\right) ) =\sup_{   x < - \kappa_n}  \Phi \left( x\right)   =  \Phi \left( - \kappa_n\right)   \leq  C_{\rho, 4} \Big(  \delta_m +\gamma_m |\ln \gamma_m| +\varepsilon_m^{ \rho(2-\rho)/4 }+\sqrt{\frac{m}n}   \Big ).
\end{eqnarray}
Combining the inequalities (\ref{ineq010})-(\ref{ineq050}) together, we obtain the desired inequality.

\subsection{Proof of Proposition \ref{pro3.3}  }\label{sec4.5}
 We only need to show that the quantities $\gamma_m$ and $\delta_m$ can  be dominated   via the quantities
   $\eta_{1, n}$ and $ \eta_{2, n} .$
By the definition of $\gamma_m$, it is easy to see that
   \begin{eqnarray*}
 \gamma_m &\leq&  \frac{1}{m^{1/2} \sigma  } \sum_{j=1}^{\infty} \frac{1}{ j^{3/2}}  \Big ( \sum_{i=1}^{mj} \eta_{1, i} \Big ).
  \end{eqnarray*}
  Thus, when $ \eta_{1, n}  =O(  n^{-\beta})$ for some $\beta>1$, it holds
  $$\gamma_m  =  O(1/ m^{1/2}).$$
  Next, we give an estimation for  $\delta_m$.
 It is obvious that
 \begin{eqnarray*}
 \|\mathbb{E}[ S_m |\mathcal{F}_0] \|_\infty \leq  \sum_{i=1}^{m} \eta_{1, i}
\end{eqnarray*}
and
 \begin{eqnarray*}
 \Big \| \frac1{m\sigma^2}  \mathbb{E}[ S_m^2| \mathcal{F}_0]- 1\Big \|_\infty  \leq  \frac{1}{m  \sigma_n^2  }\Big(\| \mathbb{E}[ S_m^2| \mathcal{F}_0] - \mathbb{E}[ S_m^2 ] \|_\infty  +| \mathbb{E}[ S_m^2 ]   -m\sigma^2 |\Big) .
\end{eqnarray*}
Clearly, it holds
\begin{eqnarray*}
 &&  \|  \mathbb{E}[ S_m^2| \mathcal{F}_0] - \mathbb{E}[ S_m^2 ]  \|_\infty  \leq  \sum_{i=1}^{m}\| \mathbb{E}[X_i^2 |\mathcal{F}_0] -\mathbb{E}[X_i^2 ] \|_\infty  \\
 && \ \ \ \ \ \ \ \ \ \ \ \ \ \ \ \ \ \  \ \ \ \ \ \ \ \ \   \ \ \ \ \ \     +\ 2 \sum_{i=1}^{m-1}\sum_{j=i+1}^{m}\| \mathbb{E}[X_iX_j |\mathcal{F}_0] -\mathbb{E}[X_iX_j ] \|_\infty . \\
\end{eqnarray*}
 Splitting the last sum as
   \begin{eqnarray*}
  \sum_{1\leq i \leq m/2}    \sum_{  i+1 \leq j \leq 2i }    +  \sum_{1\leq i \leq m/2}    \sum_{  2i+1 \leq j \leq m  }   +   \sum_{ m/2 \leq i \leq m-1}    \sum_{  i+1 \leq j \leq m }  ,
  \end{eqnarray*}
 by the condition $\max_{i=1,2 }\{   \eta_{i, n}  \} =O(  n^{-\beta}),$
  we infer that
 \begin{eqnarray*}
 \|  \mathbb{E}[ S_m^2| \mathcal{F}_0] - \mathbb{E}[ S_m^2 ]  \|_\infty  \leq  C_1\bigg( \sum_{i=1}^{m } i^{-\beta} +  \sum_{1\leq i \leq m/2}  i\eta_{2, i} + \|X_0\|_\infty\!\!\sum_{1\leq i \leq m/2}    \sum_{  j \geq  i  }\eta_{1, j} +   m \!\! \sum_{  i \geq m/2  }  \eta_{2, i} \bigg),
  \end{eqnarray*}
  Notice that
   $$|  \mathbb{E}[X_iX_0 ] | = |  \mathbb{E}[ X_0 \mathbb{E}[X_i|\mathcal{F}_{0} ] ] | \leq \| X_0 \|_\infty \| \mathbb{E}[X_i|\mathcal{F}_{0} ] \|_\infty.$$
By $ \eta_{1, n}  =O(  n^{-\beta}), \beta> 1$,   it is easy to see that
\begin{eqnarray*}
 | \mathbb{E}[ S_m^2 ]   -m\sigma^2 | &\leq&  \sum_{i=1}^m \Big|\sum_{j=1}^m \mathbb{E}[X_iX_j ] - \sum_{j=-\infty}^\infty \mathbb{E}[X_iX_j ] \Big| \\
&=& \sum_{i=1}^m \Big|\! \sum_{j=-\infty}^0\mathbb{E}[X_iX_j ]  +  \sum_{j=m+1}^{ \infty} \mathbb{E}[X_iX_j ]  \Big| \\
&\leq& \|X_0\|_\infty \sum_{i=1}^m \Big(\sum_{j=-\infty}^{-i}  O(  |j|^{-\beta} )   +  \sum_{j=m+1-i}^{ \infty} O(  j^{-\beta} )  \Big)   \\
  & \leq&   C_2 \sum_{i=1}^{m }   i^{-\beta}  \\
& \leq& C_3 \,.
\end{eqnarray*}
Hence, it holds
 \begin{eqnarray*}
 \delta_m^2 \leq  \frac{C_1}{m  \sigma_n^2  }  \Bigg[ \Big ( \sum_{i=1}^{m} \eta_{1, i} \Big )^2 +\! \sum_{i=1}^{m } i^{-\beta}+\!\!\!\sum_{1\leq i \leq m/2}  i\eta_{2, i} + \|X_0\|_\infty\!\!\sum_{1\leq i \leq m/2}    \sum_{  j \geq  i  }\eta_{1, j} +   m \!\! \sum_{  i \geq m/2  }  \eta_{2, i} +C_4 \Bigg].
  \end{eqnarray*}
Then, taking into account that $ \max_{i=1,2}\{   \eta_{i, n}  \} =O(  n^{-\beta}),$  we have
  \begin{displaymath}
 \delta_m =  \left\{ \begin{array}{ll}
  O( m^{-(\beta-1)/2} ), & \textrm{\ \ \  if $\beta \in (1, 2)$,}\\
 O( m^{-1/2} \sqrt{\ln m} ),  & \textrm{\ \ \ if $\beta = 2$,} \\
O( m^{-1/2} ), & \textrm{\ \ \  if $\beta > 2$.}
\end{array} \right.
\end{displaymath}
 By point 2 of Remark \ref{re01},
we have  $\varepsilon_m= O(  m /n^{1/2})$.
 If 
 $\beta\geq 3/2,$  then equality  (\ref{thls}) with $m=\lfloor n^{2/7}\rfloor$  holds uniformly  for $ 0 \leq x = o(n^{1/14} /\sqrt{\ln n} )$ as $n \rightarrow \infty.$
 If 
 $\beta \in (1, 3/2),$  then equality  (\ref{thls}) with $m=\lfloor n^{1/(3\beta-1)}\rfloor$ holds uniformly  for $ 0 \leq x = o(n^{(\beta-1)/(6\beta-2)}   )$ as $n \rightarrow \infty.$ This completes the proof of points [i] and [ii].

 To prove [iii],   notice   that  $m :=m(n)\rightarrow \infty$ and $n^{ 1/2} /m \rightarrow \infty$  imply   $\varepsilon_m, \gamma_m, \delta_m \rightarrow 0 $ as $n \rightarrow \infty.$
Then, point [iii] follows by Corollary  \ref{co0}.
 \hfill\qed

\section*{Acknowledgements}
The authors are deeply indebted to the editor and the anonymous referee  for their helpful comments. The work
has been supported by  the National Natural Science Foundation of China (Grant nos.\ 11601375,  11971063,   11571052 and  11731012).
The work has also benefited from the support of the Centre Henri Lebesgue (CHL, ANR-11-LABX-0020-01).


\end{document}